\documentclass[11pt]{article}

\usepackage[top=1.1in, bottom=.9in, left=1in, right=0.98in]{geometry}

\usepackage{amsmath,amssymb,amsthm}
\usepackage{graphicx}
\usepackage{enumerate}
\usepackage{bm}
\usepackage{xcolor}
\usepackage{mdframed}
\usepackage{bbm}
\usepackage[space]{cite}

\usepackage{tocloft}
\usepackage{hyperref}
\hypersetup{
	colorlinks = true,
	linkcolor = blue,
	citecolor = blue,
	urlcolor = blue,
	filecolor = black,
	linktoc = all,
}

\newtheoremstyle{style1} %
  {0pt} 
  {0pt} 
  {}    
  {0pt} 
  {\bfseries} 
  {.} 
  {1em} 
  {\thmname{#1}\thmnumber{ #2}\thmnote{ (#3)}} 

\definecolor{proofline}{gray}{0.7}
  
\newtheoremstyle{proofstyle} %
  {8pt} 
  {0pt} 
  {\small} 
  {0pt} 
  {\itshape} 
  {.} 
  {1em} 
  {{\color{proofline} \rule{6in}{.2pt} } \\[.6pt] Proof} 

\theoremstyle{style1}
\newtheorem{defn}{Definition}
\newtheorem{thm}[defn]{Theorem}
\newtheorem{lem}[defn]{Lemma}
\newtheorem{cor}[defn]{Corollary}
\newtheorem{prop}[defn]{Proposition}

\theoremstyle{proofstyle}
\newtheorem{pf}{Proof}

\definecolor{lg}{gray}{0.93}

\newmdenv[
    linewidth=0pt,
    backgroundcolor=lg,
    innerleftmargin=5pt,
    innerrightmargin=3pt,
    innertopmargin=3pt,
    innerbottommargin=4pt,
]{gbox}

\newcommand{\e}{e}
\newcommand{\E}{E}
\newcommand{\dd}{\delta}

\newcommand{\cd}{\cdot}
\newcommand{\w}{\wedge}
\newcommand{\grade}[1]{\langle #1 \rangle}
\newcommand{\pp}{\partial}
\newcommand{\ph}{\varphi}
\newcommand{\G}{\Gamma}
\newcommand{\lb}[2]{[#1,#2]}
\newcommand{\wcd}{\hspace{2.14pt} \cd \hspace{2.14pt}}
\newcommand{\myspace}{6pt}

\newcommand{\Gbar}{\bar{\Gamma}}
\newcommand{\GA}{\mathcal{A}}
\newcommand{\D}{\bm{D}}
\newcommand{\del}{\nabla}
\newcommand{\f}{\bm{f}}
\newcommand{\tr}{\textrm{tr}}
\newcommand{\rot}{\textrm{rot}}
\newcommand{\id}{\mathbbm{1}}
\newcommand{\Dbar}{\bar{D}}
\newcommand{\Dres}{\tilde{D}}
\newcommand{\DD}{\del}
\newcommand{\DDD}{\bm{\nabla}}
\newcommand{\dhat}{\hat{d}}
\newcommand{\wh}{\widehat}

\usepackage{scalerel,stackengine}
\stackMath
\newcommand\reallywidehat[1]{%
\savestack{\tmpbox}{\stretchto{%
  \scaleto{%
    \scalerel*[\widthof{\ensuremath{#1}}]{\kern-.6pt\bigwedge\kern-.6pt}%
    {\rule[-\textheight/2]{1ex}{\textheight}}
  }{\textheight}%
}{0.5ex}}%
\stackon[1pt]{#1}{\tmpbox}%
}

\begin{document}

\begin{center}
{\large \textbf{Geometric calculus on pseudo-Riemannian manifolds}} \\
\vspace*{10pt} Joseph C. Schindler \\
\vspace*{6pt} \textit{\small SCIPP and Department of Physics, University of California Santa Cruz, Santa Cruz, CA, USA}\\
\vspace*{0pt} {\small Email: \texttt{jcschind@ucsc.edu} }\\
\vspace*{6pt} {\small September 2020 } \\
\end{center}



\begin{abstract}
This article provides a pedagogically oriented introduction to geometric (Clifford) calculus on pseudo-Riemannian manifolds. Unlike usual approaches to the topic, which rely on embedding the geometric algebra either within a tensor algebra or within a vector manifold framework, here we define geometric calculus directly, by elementary methods. In particular we use an axiomatic approach that directly parallels textbook introductions to general relativity and pseudo-Riemannian geometry, so that no structure outside the metrical Clifford bundle of the manifold need be introduced. On this basis we develop the full theory of differential calculus for vector, multivector, and tensor fields.
\end{abstract}


\centerline{
\parbox{0.9\textwidth}{ \setlength{\parindent}{0pt} 
}}

\vspace{8pt}

\centerline{
\parbox{0.9\textwidth}{ \setlength{\parindent}{0pt} 
}}

{
\small
\tableofcontents
}


\clearpage

\section{Introduction}

Methods of geometric (Clifford) algebra have useful applications throughout physics and mathematics~\cite{Clifford1878,hestenes87,doran07,macdonald2017}. Applying these methods to extend traditional vector and tensor calculus is the domain of geometric calculus. Like tensor calculus, it is natural to study geometric calculus in the context of pseudo-Riemannian manifolds (that is, smooth manifolds with a nondegenerate metric of arbitrary signature). 

Typically geometric calculus is applied to this problem using one of two approaches. One, the algebraic approach, is typified by textbooks such as \cite{lawson1989,crumeyrolle1990,rodrigues2007}, where a Clifford bundle is defined as a quotient of the tensor bundle of the manifold. The other, the vector manifold approach, is typified by textbooks such as~\mbox{\cite{hestenes1984,doran2003}}, where the manifold is itself taken to be embedded within a geometric algebra. In either case, some additional structure beyond the manifold and its Clifford bundle must be used to define the formalism.

Here we will take a slightly different approach where no additional structure is necessary, in order to more closely parallel standard textbook treatments of pseudo-Riemannian geometry and general relativity. The approach is based on a minimal set of axioms for a multivector directional derivative (MDD) operator (that is, an extension of the affine connection to act on multivector fields), encapsulated in Definition~\ref{thm:def:directional-deriv}. The entirety of the framework is derived from these axioms, just as standard Riemannian geometry is based on the axioms for an affine connection~\cite{lee02,lee97}.

Within this approach, an important point is that it is not trivial to show that the MDD axioms are self-consistent---that is, it is not obvious that any operator satisfying those axioms exists. Therefore, the direct constructive proof of Theorem~\ref{thm:thm:mdd-exists} (existence) is a key result underpinning the formalism. Once this existence is established, many any other useful properties can be proved from the axioms, including that MDDs are grade-preserving, follow product rules for the geometric, dot, and wedge products, and can be used to define the tensor covariant derivative, differential forms exterior derivative, and multivector gradient.

This framework has some useful features. It naturally incorporates arbitrary bases (which may be holonomic, orthonormal, or neither), metrics of arbitrary signature, and connections both with and without torsion. MDDs will be seen to have powerful properties, and the tensor derivative will be derived from MDDs based on a chain rule for linear functions. Then straightforward generalizations of the vector calculus gradient, divergence, and curl are introduced, and it is shown that the theory of differential forms is subsumed.

An important aspect to notice about the geometric algebra approach is that it works only in terms of tangent vectors and tangent multivectors, but dual vectors (\textit{i.e.}~covectors) will never be introduced. Given a vector basis $e_i$, rather than defining a dual basis in the dual space, one introduces a ``reciprocal'' basis $e^i$ (obeying $e^i \cd e_j = \dd^i_j$) of vectors \textit{in the same space}. This is possible because, given the metric, there is a canonical isomorphism identifying covectors with vectors. This simplifies certain parts of the theory, but also means these methods can be used only when a metric is present. Nonetheless, it is possible to analyze smooth manifolds by assuming an arbitrary metric and forming only metric-independent results within the framework.

The structure of the paper is as follows.

Geometric algebra (the algebra of ``multivectors'') is reviewed in Appendix~\ref{sec:GA}, while Section~\ref{sec:survey} provides a brief survey of the framework that will be found in the main text, as well as some further motivation for the framework.

Sections~\ref{sec:notation}--\ref{sec:scalar-vector} establish conventions, define a geometric tangent space to the manifold, and review the familiar theories of scalar fields, tangent vector fields, affine connection, and Lie bracket. The purpose of this review is both pedagogical, and to provide conventions for the later sections. In the interest of a self-contained treatment, detailed proofs are given even in the review.

The main results occur in Sections~\ref{sec:multivector}--\ref{sec:tensor-fields}, where we establish the differential theories of multivectors, tensors, and differential forms, all in terms of the MDD axioms. Various details and several longer proofs appear in the appendix.

\section{Survey and motivation}
\label{sec:survey}

Geometric algebra (the algebra of ``multivectors'') is reviewed in Appendix~\ref{sec:GA}. 

The following is a summary of the basic formalism to be established for geometric calculus on pseudo-Riemannian manifolds:
\begin{itemize}
    \item Let $\e_i$ be an arbitrary set of basis vector fields (i.e. an arbitrary frame). There exists a reciprocal frame $\e^i$ defined by $\e^i \cd \e_j = \dd^i_j$. Metric coefficients and Lie bracket coefficients in the basis are defined by
    \begin{equation*}
        \e_i \cd \e_j = g_{ij},
        \qquad \qquad
        [\e_i, \e_j ] = L_{ijk} \, \e^k \, ,
    \end{equation*}
    and connection coefficients in the basis are defined by
    \begin{equation*}
        D_{\e_i} \e_j = \Gamma_{ijk} \, \e^k
    \end{equation*}
    for an affine connection $D$.
    \item There are three important types of special basis: orthonormal, holonomic, and gradient. An orthonormal basis is one such that $g_{ij}= \eta(i) \, \dd_{ij}$ with $\eta(i)=\pm 1$. A holonomic basis is one such that $L_{ijk}=0$. A gradient basis is one such that $\e_i = \DDD \ph_i$  (gradient $\DDD$ defined later) for some set  $\ph_i$ of scalar fields. An orthonormal basis is reciprocal to an orthonormal basis. A holonomic basis is reciprocal to a gradient basis, and vice-versa.
    \item We use a special notation for coordinate bases. For a coordinate system $x^i$, coordinate bases are denoted $e(x_i)$ (this would usually be denoted $\frac{\pp}{\pp x^i}$), and are holonomic. Coordinate derivatives of a scalar field are then written $\pp_{e(x_i)} \ph= \frac{\pp \ph}{\pp x^i}$. The reciprocal basis to $e(x_i)$ is the coordinate gradient basis $dx^i = \DDD x^i = g^{ij} e(x_i)$ obtained from the gradient of the scalar coordinate functions. These are both vector bases, no dual space was introduced. 
    \item In an arbitrary basis~$e_i$ (which might be holonomic, orthonormal, or neither), metric-compatibility and torsion-freeness of the connection are expressed by
    \begin{equation*}
        \Gamma_{ijk} + \Gamma_{ikj} = \pp_{\e_i} g_{jk},
        \qquad \qquad
        \Gamma_{ijk} - \Gamma_{jik} = L_{ijk}.
    \end{equation*}
    Thus if general connection coefficients are written as
    \begin{equation*}
    \Gamma_{ijk}  
        = \tfrac{1}{2} \left( \pp_{\e_i}g_{jk} - \pp_{\e_k}g_{ij} + \pp_{\e_j}g_{ki} \right)
        + \tfrac{1}{2} \left( L_{ijk} - L_{jki} + L_{kij} \right)
        + \chi_{ijk}
    \end{equation*}
    it immediately follows that the Levi-Civita connection is given by $\chi_{ijk}=0$. In general the $\chi_{ijk}$ are called ``contorsion coefficients''.
    \item The affine connection is extended to a ``multivector directional derivative'' (MDD) by demanding the product rule
    \begin{equation*}
        D_{a} (AB) = (D_{a} A) \, B + A \, (D_{a} B)
    \end{equation*}
    hold on all multivectors $A,B$. MDDs will be defined axiomatically, and their existence proved by construction. It is also proved from the axioms that MDDs are grade-preserving, and that a similar product rule also holds for the dot and wedge products.
    \item MDDs are extended to act on tensor fields by a direct chain-rule computation. For example if $T(A,B)$ is a multilinear function (tensor) with multivector inputs $A,B$ then, by the chain rule, its tensor derivative $D_a T$ is defined by
    \begin{equation*}
        D_a T(A,B) = D_a(T(A,B)) - T(D_a A, B) - T(A,D_a B)
    \end{equation*}
    where $a$ is a vector and derivatives on the right hand side are MDDs. The relation of this expression to the chain rule is clarified in the text. In terms of tensor coefficients, this becomes the usual
    \begin{equation*}
        D_{i} T^{jk}  = \pp_i T^{jk} + (\G_{ilm} \, g^{mj}) \, T^{lk}  + (\G_{ilm} \, g^{mk}) \, T^{jl} \; .
    \end{equation*}
    This chain rule definition of $D_a T$ is equivalent to the usual covariant derivative definition.
    \item An MDD operator is not a tensor, since in general $D_a(\ph A) \neq \ph D_a(A)$ for scalar fields $\ph$. However, the difference $(D-\tilde{D})$ between two MDDs is a tensor. Hence $\Gamma_{ijk}$ do not transform like tensor coefficients, but $\chi_{ijk}$ do.
    \item In this formalism, one normally thinks of the metric coefficients, rather than the metric tensor. But if desired, the metric tensor can be formally defined by $\hat{g}(a,b)=a \cd b$. Then its tensor derivative is $D\hat{g}=0$ by metric compatibility.
    \item The gradient operator $\D$ associated with an MDD $D$ is defined by 
    $$\D A = \e^i \, D_{\e_i} A$$ 
    on a multivector field $A$ and shown to have some useful properties, including that gradient 
    \begin{equation*}
        \D A = \D \cd A + \D \w A
    \end{equation*}
    is equal to the divergence ($\D \cd A \equiv \e^i \cd D_{\e_i} A$) plus the curl ($\D \w A \equiv \e^i \w D_{\e_i} A$). Multivector fields with zero gradient are the $n$-dimensional analogue of complex analytic functions. 
    \item The unique torsion-free MDD is defined as $\del$ and its gradient operator $\DDD = \e^i \, \del_{\e_i}$ is a special case of $\D$ above. It is proved that\\[-9pt]
    \begin{equation*}
        d = \DDD \w
    \end{equation*}
    is completely equivalent to the exterior derivative of differential forms. Consequently, $\DDD \w$~is independent of the metric even though $\DDD$ is metric-compatible. Identifying the exterior subalgebra of geometric algebra with the space of forms, the theory of differential forms is totally subsumed.
\end{itemize}

One motivation for the approach outlined above is that it can help provide intuitive descriptions of physical systems. For an example, consider three ways to express the Lagrangian density for relativistic electrodynamics~\cite{doran07}: \\[15pt]
\begin{tabular}{p{10mm}p{40mm}p{46mm}p{45mm}}
    &
    \textrm{Tensor calculus:}
    & 
    $\mathcal{L} = - \tfrac{1}{4} \, F^{\mu\nu} F_{\mu\nu} + A^{\mu} J_{\mu}$  
     & 
     $F_{\mu\nu} = \del_{\mu} A_{\nu} - \del_{\nu} A_{\mu}$  
     \\[10pt]
     &
    \textrm{Differential forms:}
    & 
    $\mathcal{L} = - \tfrac{1}{2} \, F \w * F + A \w *J$  
     & 
     $F=dA$  
     \\[10pt]
     &
    \textrm{Geometric calculus:}
    & 
    $\mathcal{L} = -\tfrac{1}{2} \, F \cd F + A \cd J$  
     & 
     $F= \DDD \w A$ .
\end{tabular} \\[15pt] 
The corresponding equations of motion (Maxwell's equations) are \\[15pt]
\begin{tabular}{p{10mm}p{41mm}p{30mm}p{45mm}}
    &
    \textrm{Tensor calculus:}
    &
    $\del_\mu F^{\mu\nu}= J^{\nu}$ 
     & 
     $\epsilon^{\alpha\beta\mu\nu} \; \del_{\alpha} F_{\mu\nu} = 0$ 
     \\[10pt]
     &
    \textrm{Differential forms:}
    &
    $*\,d*F= J$ 
     &
    $dF=0$ 
    \\[10pt]
     &
    \textrm{Geometric calculus:}
    &
    $\DDD \cd F = J$ 
     &
    $\DDD \w F = 0$ .
\end{tabular} \\[15pt]
The geometric calculus operations are just what they look like: they generalize the divergence and curl from vector calculus. The bottom equations are especially useful because they express both metric dependent (divergence term) and independent (curl term) concepts concisely, and help clarify how these expressions fit into a larger theoretical framework.

The geometric calculus form of these equations simplifies further, to the single expression
\begin{equation*}
    \DDD F = J.
\end{equation*}
The divergence and curl equations are the vector and trivector parts respectively (see later sections). This is more than a notational change---the operator $\DDD$ can be inverted by an integral formula, allowing a direct calculation of field strength $F$ from the source $J$. Geometric calculus also usefully unifies vector, tensor, and spinor calculus. In particular, spinors are a subset of multivectors, and Dirac's equation can be expressed using the same derivative operator~$\DDD$ appearing in Maxwell's equation~\cite{Hestenes2015}.

Throughout the rest of the article, the formalism surveyed above is developed in detail.

\section{Notation and conventions}
\label{sec:notation}

Various notations will be defined below---we collect and review some important ones here. 

Except where noted otherwise, our conventions generally align with Lee~\cite{lee02,lee97} for smooth and Riemannian manifolds, and with Macdonald~\cite{macdonald2017,macdonald11,macdonald12} for geometric algebra. Further review of geometric algebra conventions is in Appendix~\ref{sec:GA}. 

Symbols use will generally abide by the following conventions:

\begin{equation*}
    \begin{array}{lcl}
        D_a, \D & \quad & \textrm{directional derivative and gradient} \\
        \DD_a, \DDD & \quad & \textrm{torsion-free directional derivative and gradient} \\
        d & \quad & \textrm{torsion-free curl ($d=\DDD \w$) i.e. exterior derivative}\\
          & \quad &  \\
        \alpha, \beta, \varphi, \ldots & \quad & \textrm{scalar fields} \\
        a, b, \ldots & \quad & \textrm{vector fields} \\
        A, B, \ldots & \quad & \textrm{multivector fields} \\
        T, S, \ldots & \quad & \textrm{tensor fields} \\
        \e_i, \e^j, \E_i, \E^j, \ldots & \quad & \textrm{basis vector fields (basis frames)} \\
        \e_J, \e^K, \E_J, \E^K, \ldots & \quad & \textrm{basis multivector fields (multivector frames)} \\
          & \quad &  \\
        x_i\equiv x^i & \quad & \textrm{coordinate system} \\
        e(x_i) & \quad & \textrm{coordinate tangent basis (coordinate frame)} \\
    \end{array}
\end{equation*}

A few more notes will also be useful:
\begin{itemize}

\item \textit{(Derivative subscripts.)} The subscript $a$ in $D_a$ and $\DD_a$ is a vector, not a basis index. Rather, $D_{\e_i}$ would be the derivative in the~$\e_i$ basis direction. Several types of directional derivative appear (affine, multivector, tensor), but there is no danger of ambiguity (except as addressed in Definition \ref{thm:def:tensor-derivative-notation}), since all others derive from the multivector directional derivative, and each is made clear in context. 

\item \textit{(Index placement.)} The upper and lower indices in \textit{basis frames}, for example in $\e_i$ and $\e^j$, refer to a pair of reciprocal bases (that is, two sets of basis vectors mutually chosen such that $\e^i \cd \e_j = \dd^i_j$). However, upper and lower index placement for \textit{coordinates} is meaningless, so $x_i \equiv x^i$. The symbols $\e_i$~and~$\E_i$  refer to arbitrary basis frames;  whether such frames are orthonormal, holonomic, or neither depends on context. Likewise, the symbol $g_{ij}$ may refer to the metric in any basis, depending on context, and is not reserved for the coordinate basis metric. 

\item \textit{(Summation convention.)} We assume the Einstein summation convention, so that repeated indices are always summed unless it is explicitly stated otherwise (this excludes indices inside a function argument, for example $i$ in $\eta(i)$, so that $\eta(i) a^i \dd_{ij}$ would be summed but  $\eta(i) \dd_{ij}$ would not). 

\item \textit{(Global vs. local fields.)} Fields (scalar, vector, etc.) on a manifold will generally be assumed to be defined only locally in a neighborhood of some point, unless the field is explicitly said to be defined globally. 

\item \textit{(Smoothness.)} The term \textit{smooth} implies that at least as many derivatives exist as are needed for any calculation, but need not imply infinite differentiability. Fields are assumed smooth unless stated otherwise.

\end{itemize}

\section{Geometric tangent space}
\label{sec:GM}

Let $M$ be a smooth manifold with a nondegenerate bilinear metric of arbitrary signature (\textit{i.e.}~a pseudo-Riemannian manifold). The tangent space $T_p  M$ (constructed in the usual way~\cite{lee02}) can be extended to a \textit{geometric tangent space} $GT_p M$ at each point, using the pseudo-Riemannian metric as the geometric algebra dot product. Elements of $GT_p M$ are \textit{tangent multivectors}, which contain tangent vectors in $T_p M$ as a subset. $MV\!F(M)$ denotes the space of smooth multivector fields on~$M$ (either local or global depending on context). The smooth structure of $GT_p M$ is discussed in Appendix~\ref{sec:smooth}.

We adopt some nonstandard notational conventions for tangent vectors and bases, both to reduce the emphasis usually placed on coordinate systems, and to better incorporate tangent vectors into the larger multivector framework. We mainly work in terms of completely arbitrary vector bases $e_i$ (which may be orthonormal, holonomic, or neither) with the notation:\\[-6pt]
\begin{equation*}
    \begin{array}{rcccccc}
        & \quad &
        \textrm{\small Arbitrary Basis}
        & \quad &
        \textrm{\small Tangent Vector}
        & \quad &
        \textrm{\small Directional Derivative}
        \\[4pt]
        \textrm{Standard Notation:}
        & \quad &
        \textit{\small(uses coordinates)}
        & \quad &
        X 
        & \quad &
        X \ph 
        \\[5pt]
        \textrm{New Notation:}
        & \quad &
        \e_i
        & \quad &
        a
        & \quad &
        \pp_a \, \ph
    \end{array}
\end{equation*}\\[2pt]
To distinguish coordinate bases then requires a special notation (with coordinates $x^i$):\\[-4pt]
\begin{equation*}
    \begin{array}{rcccccc}
        & \quad &
        \textrm{\small Coordinate Basis}
        & \quad &
        \textrm{\small Tangent Vector}
        & \quad &
        \textrm{\small Directional Derivative}
        \\[4pt]
        \textrm{Standard Notation:}
        & \quad &
        \frac{\pp}{\pp x^i}
        & \quad &
        X = X^i \; \frac{\pp}{\pp x^i}
        & \quad &
        X \ph = X^i \; \frac{\pp \ph}{\pp x^i} 
        \\[8pt]
        \textrm{New Notation:}
        & \quad &
        \e(x_i)
        & \quad &
        a= a^i \; \e(x_i)
        & \quad &
        \pp_a \, \ph  = a^i \; \pp_{\e(x_i)} \ph = a^i \; \frac{\pp \ph}{\pp x^i}
    \end{array}
\end{equation*}\\[0pt]
Note that for coordinates (but not basis vectors) index placement is meaningless ($x_i \equiv x^i$) and is chosen for notational convenience. The above notation involves the identification
\begin{equation*}
    \pp_{e(x_i)} = \tfrac{\pp}{\pp x^i}
\end{equation*}
so that $e(x_i)$ is the tangent vector pointing in the direction of the coordinate partial derivative. More broadly we write
\begin{equation*}
    \pp_a \ph
\end{equation*}
rather than $a(\ph)$ for the directional derivative of scalar field $\ph$, emphasizing that $T_p M$ is the space of ``directional derivative directions''.

We denote vector ($T_p M$) basis frames by $\e_i$ (lowercase latin index). From vector basis frames can be constructed multivector ($GT_p M$) basis frames $e_J$ (uppercase latin index), defined by
\begin{equation*}
    \e_J=\e_{j_0} \w \ldots \w \e_{j_n}
\end{equation*}
as in Appendix~\ref{sec:smooth}. Arbitrary vector fields $a = a^i \e_i$ and multivector fields $A = A^J \e_J$ can be defined in terms of a frame. The definition of smoothness for fields and frames is given in Appendix~\ref{sec:smooth}. All fields and frames are assumed smooth unless stated otherwise.

\section{Basis frames}
\label{sec:frames}

A basis frame (or simply ``basis'', when the fact we are discussing fields is implied, as below) is a sufficiently smooth (see Appendix~\ref{sec:smooth}) set of vector fields that forms a vector basis at each point. 

An arbitrary basis $\e_i$ is characterized by two important quantities, the metric coefficients~$g_{ij}$ and Lie bracket coefficients~$L_{ijk}$ (Lie bracket defined later) defined by
\begin{equation*}
    g_{ij} = \e_i \cd \e_j \, ,
    \qquad \qquad \qquad
    [\e_i, \e_j] = L_{ijk} \, \e^k .
\end{equation*}
The reciprocal basis $\e^i$ is defined by $\e^i \cd \e_j = \dd^i_j$ as always (see Appendix~\ref{sec:GA}), so that $\e^i = g^{ij} \e_j$ expresses the reciprocal basis using the matrix inverse of the metric. In terms of the basis and reciprocal basis one can expand an arbitrary vector field $a$ by
\begin{equation*}
    a = a^i \, \e_i = a_i \, \e^i
\end{equation*}
where $a^i = a \cd \e^i = g^{ij} \, a_j$ and $a_i = a \cd \e_i = g_{ij} \, a^j$.

There are three important types of special basis: orthonormal, holonomic, and gradient bases. They are defined as follows:\\[-16pt]
\begin{center}
\begin{tabular}{rcccccccc}
    & \quad &
    Arbitrary
    & \quad &
    Orthonormal
    & \qquad & 
    Holonomic
    & \qquad & 
    Gradient
    \\[4pt]
    Definition:
    & \quad &
    $\e_i$
    & \quad &
    $g_{ij} = \eta(i) \, \dd_{ij}$
    & \qquad & 
    $L_{ijk} = 0$
    & \qquad & 
    $\e_i = \DDD \ph_i$
    \\[6pt]
    Reciprocal to:
    & \quad &
    $\e^i = g^{ij} \, \e_j$
    & \quad &
    \textit{\small Orthonormal}
    & \qquad & 
    \textit{\small Gradient}
    & \qquad & 
    \textit{\small Holonomic}
\end{tabular}
\end{center}

\noindent
To clarify the notation above, a basis is a gradient basis if there is a set of scalar fields $\ph_i$ such that each basis field $\e_i$ is the gradient of one of the $\ph$. The gradient operator $\DDD$ will be defined later. For orthonormal bases the function $\eta(i)=\pm 1$ (for each index) determines the metric signature, and often one writes the orthonormal metric as $\eta_{ij} \equiv \eta(i) \, \dd_{ij}$. It is shown below that every holonomic basis is reciprocal to a gradient basis, and vice versa. Moreover, it is easy to verify that the reciprocal to an orthonormal basis is also orthonormal. In particular an orthonormal basis obeys $\e^i = \eta(i) \, \e_i = \pm \e_i$, so an orthonormal basis in Euclidean signature is self-reciprocal.

Given a coordinate system $x^i$, the coordinate basis $\e(x_i)$ is holonomic. On the other hand, each individual coordinate $x^i$ can be treated as a scalar field. The gradient works out to $\DDD x^i = g^{ij} \, e(x_j)$, ensuring that $\DDD x^i$ forms a basis reciprocal to $e(x_i)$. It will be shown later on that for scalar fields $\DDD \ph = \DDD \w \ph = d \ph$, so we typically write the coordinate gradient in the more suggestive notation~$dx^i$. Thus the situation for a coordinate system $x^i$ can be summarized by
\begin{center}
\begin{tabular}{rcccc}
    & \quad &
    Coordinate Basis
    & \qquad \qquad & 
    Coordinate Gradient Basis
    \\[4pt]
    & \quad &
    $\e(x_i)$
    & \qquad & 
    $dx^i = \DDD x^i = g^{ij} \, e(x_j) $
    \\[4pt]
    & \quad &
    \textit{\small Type: Holonomic}
    & \qquad & 
    \textit{\small Type: Gradient}
\end{tabular}
\end{center}
which are mutually reciprocal with the reciprocality relation
\begin{equation*}
    e(x_i) \cd dx^j = \dd^j_i \, .
\end{equation*}
Note that the $dx^i$ basis here is \textit{still a vector basis} living in the same vector space as the coordinate basis $e(x_i)$. No dual vector space was introduced. The coordinate gradient basis has the usual properties of the dual basis, but with a simpler interpretation: the vector field $dx^i \equiv \DDD x^i$ is precisely the gradient of the scalar coordinate functions. Given any two separate coordinate systems $x^i,y^j$ it follows from the definitions that
\begin{equation*}
    e(x_i) \cd dy^j = \tfrac{\pp y^j}{\pp x^i}
\end{equation*}
leading to the usual formulas for change of coordinates.

Every holonomic basis is equivalent to a coordinate basis for some coordinate system~\cite{Schutz1980}. Similarly, the scalar fields defining a gradient basis always provide a coordinate system (since the Jacobian is constrained by the basis condition). Therefore holonomic and gradient bases come in pairs, and may always be thought of as arising from coordinate systems. This guarantees that at every point in $M$ one can choose a smooth local holonomic or gradient basis as desired. It is also straightforward to show there exists a smooth local orthonormal basis at every point.

\section{Scalar and vector fields}
\label{sec:scalar-vector}

Differential geometry begins with the directional derivative of scalar fields. The notation mentioned earlier is now formally defined.

\begin{gbox}
\begin{defn}[Scalar-field directional derivative]
\label{thm:def:scalar-deriv}
The \textit{scalar-field directional derivative} $\pp_a \ph$ of a scalar field $\ph$ in the direction of a vector $a$ at point $p$ is defined by
\begin{equation}
\pp_a \ph = a(\ph),
\end{equation}
where $a(\ph)$ is the usual action of the tangent vector $a$ on $\ph$ from the theory of smooth manifolds. Thus in coordinates $x^i$ with associated basis $e(x_i)$, the scalar-field directional derivative for a vector $a=a^i \, \e(x_i)$ is given by
\begin{equation*}
    \pp_a \ph = a^i \, \pp_{e(x_i)} \ph = a^i \, \frac{\pp \ph}{\pp x^i} \, .
\end{equation*}
\end{defn}
\end{gbox}

By definition the directional derivative is defined at a point, and returns a scalar. When a vector field is provided as the direction argument, the derivative is evaluated pointwise, and the resulting output is a scalar field. If $\alpha,\beta,\ph,\ph_1,\ph_2$ are smooth scalar fields and $a,b$ smooth vector fields, the scalar-field directional derivative has the properties 
$\pp_{(\alpha a + \beta b)} \ph
= \alpha \pp_a \ph + \beta \pp_b \ph $
and
$\pp_a (\ph_1 \ph_2) = 
(\pp_a \ph_1) \ph_2 + \ph_1 ( \pp_a \ph_2)  $. By definition, two smooth vector fields $a,b$ are equal if and only if $\pp_a \ph = \pp_b \ph$ for all $\ph$. This directional derivative notation provides a nice coordinate-independent notation to define the Lie bracket.

\begin{gbox}
\begin{defn}[Lie bracket]
\label{thm:def:lie-bracket}
For any two smooth vector fields $a$ and $b$, there exists a unique smooth vector field $c=\lb{a}{b}$, called the \textit{Lie bracket} of $a$ and $b$, such that
\begin{equation}
\label{eqn:def-lie-bracket}
\pp_{\, \lb{a}{b}} \ph = 
\pp_{   a}    \pp_{\, b} \, \ph - 
\pp_{\, b} \, \pp_{   a} \, \ph
\end{equation}
for all scalar fields $\ph$.
\end{defn}
\end{gbox}

Uniqueness follows directly from the definition of vector fields, and existence can be proven by a straightforward calculation in coordinates. That calculation reveals that in any coordinate system~$x^i$, the Lie bracket of vector fields $a=a^i \e(x_i)$ and $b=b^i \e(x_i)$ is
\begin{equation}
\lb{a}{b} = \big(a^l \, \tfrac{\pp b^k}{\pp x^l} - b^l \, \tfrac{\pp a^k}{\pp x^l} \big) \, \e(x_k).
\end{equation}
The Lie bracket has some useful properties.

\renewcommand{\myspace}{4pt}
\begin{gbox}
\begin{prop}[Lie bracket properties]
Let $a,b,c$ be smooth vector fields, and $\alpha,\beta$ be smooth scalar fields. Then \\[\myspace]
\begin{tabular}{rlr}
(i)
&
$[a+b,c] = [a,c] + [b,c]$,
&
(linearity)
\\[\myspace]
(ii)
&
$[a,b+c] = [a,b] + [a,c]$,
&
(linearity)
\\[\myspace]
(iii)
&
$[a,b] = -[b,a]$,
&
(antisymmetry)
\\[\myspace]
(iv)
&
$ [a,[b,c]] +  [c,[a,b]] + [b,[c,a]] = 0$,
&
(Jacobi identity)
\\[\myspace]
(v)
&
$ [\alpha a, \beta b] = 
\alpha (\pp_a \beta) b - \beta (\pp_b \alpha) a + \alpha \beta [a,b] $.
&
\hspace{15mm}
(scalar-field multiplier formula)
\end{tabular}
\end{prop}
\begin{pf}
The proofs rely on the linearity of and product rule for scalar-field directional derivatives. \\
(i) $\pp_{[a+b,c]} \ph = 
\pp_{a} \pp_{c} \ph + \pp_{b} \pp_{c} \ph - 
\pp_{c} \pp_{a} \ph - \pp_{c} \pp_{b} \ph =
\pp_{[a,c]+[b,c]} \ph $. \\
(ii) Likewise.\\
(iii) $\pp_{[a,b]}\ph = - (\pp_b \pp_a \ph - \pp_a \pp_b \ph) = \pp_{-[b,a]} \ph$. \\
(iv) Direct expansion of 
$\pp_{[a,[b,c]] +  [c,[a,b]] + [b,[c,a]]} \ph$ 
yields zero. \\
(v) 
$ 
\pp_{[\alpha a, \beta b]} \ph 
=
\pp_{\alpha a} \pp_{\beta b} \ph -
\pp_{\beta b} \pp_{\alpha a} \ph 
= 
\alpha \pp_a ( \beta \pp_b \ph) -
\beta \pp_{b} ( \alpha \pp_{a} \ph) 
= 
\alpha (\pp_a  \beta) \pp_b \ph +
\alpha \beta  \pp_a \pp_b \ph - 
\beta (\pp_{b}  \alpha) \pp_{a} \ph - 
\alpha \beta \pp_{b}  \pp_{a} \ph 
=
\alpha (\pp_a  \beta) \pp_b \ph -
\beta (\pp_{b}  \alpha) \pp_{a} \ph +
\alpha \beta \pp_{[a,b]} \ph
=
\pp_{\alpha (\pp_a  \beta) b - \beta (\pp_{b}  \alpha)a +
\alpha \beta [a,b]} \ph.
$
\end{pf}
\end{gbox}

Any coordinate basis has the commutators $[\e(x_i),\e(x_j)]=0$. Other basis fields, however, may have nonzero commutators. For arbitrary basis fields $\e_i$, define the \textit{commutator coefficients} $L_{ijk}$ by
\begin{equation}
[\e_i,\e_j] = L_{ijk} \e^k .
\end{equation}
It follows from antisymmetry that
\begin{equation}
L_{ijk}+L_{jik}=0
\end{equation}
for all $i,j,k$, and from the Jacobi identity that 
\begin{equation}
\sum\nolimits_{pq} g^{pq} \left(
L_{jkp} L_{qim} +
L_{ijp} L_{qkm} +
L_{kip} L_{qjm}
\right)
=
\pp_{\e_i} L_{jkm} + 
\pp_{\e_k} L_{ijm} +
\pp_{\e_j} L_{kim}
\end{equation}
for all $i,j,k,m$,  where $g_{ij} = \e_i \cd \e_j$.
Expanding in terms of the arbitrary basis fields, a direct  calculation shows that the Lie bracket of any two smooth vector fields $a=a^i \e_i$ and $b=b^i \e_i$ is
\begin{equation}
\label{eqn:arb-commutator-coeffs}
[a,b] =  
( \pp_{a} \, b^k - \pp_{\,b} \, a^k 
+ a^i b^j L_{ijm} \, g^{mk} ) 
\, \e_k .
\end{equation}
The space $V\! F(M)$ of smooth vector fields on $M$ forms an infinite-dimensional Lie algebra under the Lie bracket commutator.

A directional derivative on the space of smooth vector fields is provided by the concept of affine connection.

\begin{gbox}
\begin{defn}[Affine connection]
\label{thm:def:affine-connection}
Let 
$D: T_p M \times V\!F(M) \to T_p M$
be an operator mapping a tangent vector $a$ at $p$ and a smooth vector field $u$ in a neighborhood of $p$ to a tangent vector \mbox{$D_a u$ at $p$}.  
\mbox{$D$ is called} an \textit{affine connection} if, for all vector fields $u,v$, scalar fields $\lambda$, tangent vectors $a,b$ at $p$, and scalars $\alpha,\beta$ at $p$, it has the properties: \\[4pt]
\begin{tabular}{rlr}
(i) & $D_{(\alpha \, a + \beta \, b)} v = \alpha \, D_a v + \beta \, D_b v,  $ &
\hspace{32mm} (linearity in the direction argument)
\\[4pt]
(ii) & $D_a (u+v) = D_a u + D_a v, $ & (linearity in the field argument)
\\[4pt]
(iii) & $D_a (\lambda v) = (\pp_a \lambda) v + \lambda (D_a v). $ & (scalar-field product rule)
\end{tabular} \\[\myspace]
When $D$ is an affine connection, $D_a u$ is called the \textit{affine derivative} of $u$ in the direction $a$.
\end{defn}
\end{gbox}

By definition, the affine derivative is taken at a point (the basepoint of the direction vector), and returns a single vector at that point. However, one is free to evaluate the derivative at many points simultaneously by providing a vector field in the direction argument. Since linearity in the direction argument is evaluated pointwise, property (iii) then becomes $D_{(\alpha(x) \, a(x) + \beta(x) \, b(x))} u = \alpha(x) \, D_{a(x)} u + \beta(x) \, D_{b(x)} u  $ (where $x$ indicates a function of the coordinates), a notation which is often used in this context. It is common to insist on providing a smooth vector field as the direction argument, in which case the connection is an operator mapping 
$V\!F(M) \times V\!F(M) \to V\!F(M)$,
but the present definiton more accurately captures the role of the direction argument, and extends more directly to the concept of multivector directional derivative.

Given an arbitrary set of basis fields $\e_i$, every affine connection $D$ is determined by its \textit{connection coefficients} in that basis, defined by $\Gamma_{ijk}=(D_{\e_i} \e_j) \cd \e_k$ so that
\begin{equation}
D_{\e_i} \e_j = \Gamma_{ijk} \, \e^k .
\end{equation}
The properties of affine connections show that every set of connection coefficients defines a valid affine connection, and that two affine connections are equal if and only if their connection coefficients are equal in any and every basis.

There are two important properties used to categorize affine connections: metric-compatibility and torsion, defined as follows.

\begin{gbox}
\begin{defn}[Metric-compatibility and torsion]
\label{thm:def:metric-compatibility-torsion}
Let $D$ be an affine connection.
\\[4pt]
\begin{tabular}{r p{.9\textwidth}}
(i)
&
$D$ is called \textit{metric-compatible} if $(D_c \, a) \cd b + a \cd (D_c \, b)=\pp_c \,(a \cd b)$ for all vectors $c$ and smooth vector fields $a,b$.
\\[4pt]
(ii)
&
Let $a$ and $b$ be smooth vector fields. The \textit{torsion} $\tau(a,b)$ of $a$ on $b$ relative to $D$ is a smooth vector field  defined by
$\tau(a,b) = D_{a\,}b - D_{b\,} a - [a,b]$.
$D$ is called \textit{torsion-free} if the torsion relative to $D$ vanishes for all $a$ and $b$.
\end{tabular}
\\[2pt]
\end{defn}
\end{gbox}

Metric compatibility and torsion-freeness each correspond to a simple restriction on the connection coefficients.

\begin{gbox}
\begin{thm}[Metric-compatible and torsion-free connection coefficients]
\label{thm:thm:metric-compatible-torsion-free-gamma}
Let $D$ be an affine connection. Let $\e_i$ be an arbitrary set of basis vector fields, let $g_{ij}=\e_i \cd \e_j$ and $[\e_i,\e_j] = L_{ijk} \, \e^k$, and let 
$D_{\e_i} \e_j = \Gamma_{ijk} \,\e^k$.
\\[8pt]
\begin{tabular}{r p{.9\textwidth}}
(i)
&
$D$ is metric-compatible if and only if for all $i,j,k$ $$\Gamma_{ijk}+\Gamma_{ikj}=\pp_{\e_i} \, g_{jk} \; .$$ 
\\[-12pt]
(ii)
&
$D$ is torsion-free if and only if for all $i,j,k$ $$\Gamma_{ijk}-\Gamma_{jik} = L_{ijk} \; .$$
\\[-12pt]
\end{tabular}
\end{thm}
\begin{pf}
(i) If $D$ is metric compatible then 
$\pp_{\e_i} g_{jk} = \pp_{\e_i} (\e_j \cd \e_k) = D_{\e_i} \e_j \cd \e_k + \e_j \cd D_{\e_i} \e_k = \G_{ijk} + \G_{ikj}$. On the other hand if the formula holds then 
$
(D_c \, a) \cd b + a \cd (D_c \, b)
    = (\pp_c a^j) b^k g_{jk} + a^j (\pp_c b^k) g_{jk} + a^j b^k c^i (\G_{ijk} + \G_{ikj})
    = (\pp_c a^j) b^k g_{jk} + a^j (\pp_c b^k) g_{jk} + a^j b^k (\pp_c g_{jk})
    = \pp_c (a^j b^k g_{jk})
    = \pp_c(a \cd b)
$.
(ii) If $D$ is torsion-free then 
$D_{\e_i} \e_j - D_{\e_j} \e_i - [\e_i,\e_j] = (\G_{ijl} - \G_{jil} - L_{ijl}) \e^l = 0$. Dotting with $\e_k$ gives the result. Conversely, if the formula holds then
$
D_{a\,}b - D_{b\,} a - [a,b]
    = a^i b^j (\G_{ijl} - \G_{jil} - L_{ijl}) \e^l
    = 0
$
using (\ref{eqn:arb-commutator-coeffs}).
\end{pf}
\end{gbox}

This leads to a useful standard form for the connection coefficients.

\begin{gbox}
\begin{cor}[Standard form for connection coefficients]
\label{thm:cor:contorsion}
Let $D$ be an affine connection. Let $\e_i$ be an arbitrary set of basis vector fields, let $g_{ij}=\e_i \cd \e_j$ and $[\e_i,\e_j] = L_{ijk} \, \e^k$, and let 
$D_{\e_i} \e_j = \Gamma_{ijk} \,\e^k$. Without loss of generality, write
\begin{equation*}
\label{eqn:standard-conn-coeffs}
    \Gamma_{ijk}  
        = \tfrac{1}{2} \left( \pp_{\e_i}g_{jk} - \pp_{\e_k}g_{ij} + \pp_{\e_j}g_{ki} \right)
        + \tfrac{1}{2} \left( L_{ijk} - L_{jki} + L_{kij} \right)
        + \chi_{ijk}
\end{equation*}
where $\chi_{ijk}$ are arbitrary coefficients called the \textit{contorsion coefficients}. Then
\\[6pt]
\begin{tabular}{r p{.9\textwidth}}
(i)
&
$D$ is metric-compatible if and only if $\chi_{ijk}+\chi_{ikj}=0$ for all $i,j,k$.
\\[6pt]
(ii)
&
$D$ is torsion-free if and only if $\chi_{ijk}-\chi_{jik} = 0$ for all $i,j,k$.
\\[6pt]
(iii)
&
$D$ is both metric-compatible and torsion-free if and only if $\chi_{ijk}=0$ for all $i,j,k$.
\end{tabular}
\end{cor}
\begin{pf}
Let $A_{ijk}=\tfrac{1}{2} \left(  \pp_{\e_i}g_{jk} - \pp_{\e_k}g_{ij} + \pp_{\e_j}g_{ki}\right)$. This satisfies $A_{ijk}+A_{ikj} = \pp_{\e_i} g_{jk}$ and $A_{ijk}-A_{jik} = 0$. 
Next let $B_{ijk} = \tfrac{1}{2} \left( L_{ijk} - L_{jki} + L_{kij} \right)$. This satisfies $B_{ijk}+B_{ikj}=0$ and $B_{ijk}-B_{jik}=L_{ijk}$. 
(Note that these properties rely on $g_{ij}=g_{ji}$ and $L_{ijk}=-L_{jik}$.) 
Therefore all together one finds 
$\G_{ijk}+\G_{ikj} = \pp_{\e_i} g_{jk} + \chi_{ijk}+\chi_{ikj}$
and
$\G_{ijk}-\G_{jik} = L_{ijk} + \chi_{ijk}-\chi_{jik}$
so results (i,ii) follow immediately from Theorem~\ref{thm:thm:metric-compatible-torsion-free-gamma}.
To show (iii) first note that if $\chi_{ijk}=0$ then (i,ii) imply $D$ is metric-compatible and torsion-free. Conversely if $D$ is both metric-compatible and torsion-free then the conditions on $\chi_{ijk}$, taken together, imply $\chi_{ijk}=-\chi_{kij}$. Iterating this expression yields $\chi_{ijk}=-\chi_{kij}=\chi_{jki}=-\chi_{ijk}$ which implies $\chi_{ijk}=0$. This concludes the proof. As an aside, note that if $B_{ijk}$ were not given, it could be deduced as follows. The two desired conditions for $B_{ijk}$, applied simultaneously, give the equation $B_{ijk} = L_{ijk}-B_{jki}$. This equation can be recursively plugged into itself by substituting for $B_{jki}$ on the right hand side, and so on iteratively. Since each step cycles the indices by one slot, the process is guaranteed to eventually terminate by producing a $B_{ijk}$ on the right hand side. Indeed, after three iterations one obtains $B_{ijk} = L_{ijk} - L_{jki} + L_{kij} - B_{ijk}$ which can then be solved. An equivalent procedure can be used to deduce $A_{ijk}$. Unlike the usual derivation of the Levi-Civita connection coefficients, this method doesn't require one to guess an ungainly and unintuitive expression in order to obtain the proof.
\end{pf}
\end{gbox}

Affine connections act on vector fields. In the spirit of geometric algebra, the goal of this article is to define a similar type of operator acting on multivector fields. This is accomplished in the next section, and the resulting operator is shown to have a number of desirable and intuitive properties.

\section{Multivector fields and the multivector directional derivative}
\label{sec:multivector}

This section introduces and studies the multivector directional derivative (MDD), an operator which takes the derivative of a multivector field in the direction of a vector. A class of operators satisfying the desired axioms is shown to exist and have some useful additional properties. 

\begin{gbox}
\begin{defn}[Multivector directional derivative]
\label{thm:def:directional-deriv}
Let $A$ and $B$ be smooth multivector fields, let $a$ and $b$ be vectors based at a point $p$, and let $\alpha$ and $\beta$ be scalars \mbox{at $p$}. Let 
$D: T_p M \times MV\!F(M) \to GT_p M$
be an operator mapping a tangent vector $a$ at $p$ and a smooth multivector field $A$ in a neighborhood of $p$ to a tangent multivector \mbox{$D_a A$ at $p$}.  
\mbox{$D$ is called} a \textit{multivector directional derivative} (MDD) if it has the properties: \\[4pt]
\begin{tabular}{rlr}
(i) & $D_{(\alpha \, a + \beta \, b)} A = \alpha \, D_a A + \beta \, D_b A,  $ & (linearity in the direction argument)
\\[4pt]
(ii) & $D_a \grade{A}_0 = \pp_a \grade{A}_0 \, , $ & 
\hspace{4mm}
(scalar-field directional derivative on scalars)
\\[4pt]
(iii) & $ D_a \grade{A}_1 = \grade{D_a \grade{A}_1 }_1 \, ,$ & (preserves grade of vectors)
\\[4pt]
(iv) & $D_a (A+B) = D_a A + D_a B, $ & (linearity in the field argument)
\\[4pt]
(v) & $D_a (AB) = (D_a A) B + A(D_a B). $ & (product rule)
\end{tabular} \\[\myspace]
When $D$ is a multivector directional derivative, $D_a A$ is called the derivative of the multivector field $A$ in the direction $a$ at the point $p$.
\end{defn}
\end{gbox}

Like the affine connection and scalar-field directional derivative, the MDD is taken at a point, and returns a single multivector at that point, but can be evaluated at many points simultaneously by providing a vector field in the direction argument. If one insists on providing a smooth vector field in the direction argument, $D$ becomes an operator from \mbox{$V\!F(M) \times MV\!F(M) \to MV\!F(M)$}. 

It is not obvious from the definition whether or not an operator satisfying the above axioms exists; insisting on the product rule raises the possibility that the definition is self-contradictory or overconstrained. Fortunately, later we will see that not only does such an operator exist, but many distinct such operators exist: every metric-compatible connection can be extended to an MDD. Extending non-metric-compatible connections using the above axioms is impossible, however, since the resulting operator is self-contradictory and not well-defined.

A few important properties of directional derivatives are derived easily from the definition.

\begin{gbox}
\begin{thm}[MDD properties I]
\label{thm:thm:dd-props-1}
Let $D$ be a multivector directional derivative. Then for every vector $c$ and all smooth vector fields $a$ and $b$, the derivative $D$ obeys \\[6pt]
\begin{tabular}{rlr}
(i) & $D_{c\,} (a \wcd b) = (D_{c\,} a ) \wcd  b  +  a \wcd  (D_{c\,} b),$ & \hspace{15mm} ($D$ is \textit{metric-compatible}.)
\\[2pt]
(ii) & $D_{c\,} (a \w b) = (D_{c\,} a ) \w b  +  a \w (D_{c\,} b).$ & \hspace{15mm} ($D$ is \textit{wedge-compatible} on vectors.)
\\[2pt]
\end{tabular}
\end{thm}
\begin{pf} Expand $D_{c\,} (a \cd b) = D_{c\,} (ab + ba)/2$ with the product rule, note that $D_{c\,} a$ and $D_{c\,} b$ are vectors by definition of $D$, and regroup terms to form the right hand side. Similarly for $D_{c\,} (a \w b) = D_{c\,} (ab - ba)/2$.
\end{pf}
\end{gbox}

Thus the formalism admits only metric-compatible derivatives. This extra restrictiveness is a result of unifying the scalar and vector derivatives into a single operator such that scalar-valued products of multivectors have a well-defined derivative.

Another important property which a directional derivative may have is to preserve grade.

\begin{gbox}
\begin{prop}[Grade-preserving]
\label{thm:prop:grade-preserving}
An operator $D$ is called \textit{grade-preserving} if for every vector $a$ and smooth multivector field $A$, it obeys the equivalent conditions \\[4pt]
\begin{tabular}{rll}
(i) & $D_a \grade{A}_k = \grade{D_a \grade{A}_k}_k \qquad $ & for all grades $k$,
\\[4pt]
(ii) & $\grade{D_a A}_k = D_a \grade{A}_k \qquad $ & for all grades $k$.
\end{tabular}
\end{prop}
\begin{pf}
It is claimed that (i) is equivalent to (ii). This is shown as follows. \\
(i$\rightarrow$ii)
$\grade{D_a A}_k = 
\sum_j \grade{D_a \grade{A}_j }_k = 
\sum_j \grade{\grade{D_a \grade{A}_j }_j }_k = 
\sum_j D_a \grade{A}_j \, \dd_{jk} = 
D_a \grade{A}_k $.
\\
(ii$\rightarrow$i)
$ \grade{D_a \grade{A}_k}_k =
\grade{\grade{D_a A}_k}_k =
\grade{D_a A}_k =
D_a \grade{A}_k $.
\end{pf}
\end{gbox}

Interestingly, the product rule is strong enough to ensure that every MDD is grade-preserving.

\begin{gbox}
\begin{thm}[Grade-preserving]
\label{thm:thm:grade-pres-dd}
The axioms of Definition \ref{thm:def:directional-deriv} imply that every multivector directional derivative is grade-preserving.
\end{thm}
\begin{pf}
The proof is obtained by working in an orthonormal basis and observing a simple pattern: when evaluating the directional derivative of strings of orthonormal basis vectors, potentially non-grade-preserving terms come in pairs which together vanish by metric compatibility. Unfortunately there is no simple notation to show the proof algebraically, so the full proof is fairly involved. The full proof is given in Appendix~\ref{sec:app:mdd-grade-preserving-proof}.
\end{pf}
\end{gbox}

This result supports the intuition that a directional derivative represents a limit of differences. It has as a corollary some useful properties.

\begin{gbox}
\begin{cor}[MDD properties II]
Let $D$ be a multivector directional derivative. Then for all smooth multivector fields $A,B$ and vectors $c$, the derivative $D$ obeys \\[6pt]
\begin{tabular}{rlr}
(i)
&
$D_{c\,} (A \wcd B) = (D_{c\,} A ) \wcd  B  +  A \wcd  (D_{c\,} B),$
&
($D$ is \textit{dot-compatible})
\\[\myspace]
(ii)
&
$D_{c\,} (A \w B) = (D_{c\,} A ) \w  B  +  A \w  (D_{c\,} B),$
&
\hspace{10mm}
($D$ is \textit{wedge-compatible})
\\[\myspace]
(iii)
&
If $I$ is a smooth unit pseudoscalar field, $D_{c\,} I = 0$.
&
\end{tabular}

\end{cor}
\begin{pf} 
(i) Using linearity, the definition of the inner product, grade-preservation, and the product rule, one finds that 
$
D_c(A\cd B) 
= 
\sum_{jk} D_c \grade{\grade{A}_j \grade{B}_k}_{k-j}
=
\sum_{jk} \left( D_c \grade{A}_j \cd \grade{B}_k + \grade{A}_j \cd D_c \grade{B}_k  \right)
=
 (D_c A) \cd B + A \cd ( D_c B )
$. 
(ii) Likewise.
(iii) Let $\E_i$ be an orthonormal basis, and suppose without loss of generality that $D_{\E_i}\E_j=\gamma_{ijk} \E^k$. Metric compatibility with orthonormality implies $\gamma_{ijk}+\gamma_{ikj}=0$.  Every smooth unit pseudoscalar field is equal to $I=\pm \E_1 \w \ldots \w \E_n$, where $n$ is the dimension of $M$ so all of the $\E_i$ are represented in the product. Consider the first term in the product rule expansion for $D_{\E_i} I$, which is 
$\pm(\pm \gamma_{i1k}\E_k)\w \ldots \w \E_n$. Within this term, the $k=1$ term vanishes since $\gamma_{ijj}=0$. Meanwhile the $k\neq 1$ terms vanish by antisymmetry since every $\E_k$ for $k \neq 1$ is already represented in the wedge product. Thus the first term in the product rule expansion is zero, and all other terms in the expansion vanish for the same reason. Thus $D_a I=0$ by linearity in the direction.
\end{pf}
\end{gbox}

Having surveyed some of the properties a multivector directional derivative would have if it exists, it is time to turn to the question of existence. The issue is not trivial, but fortunately the result works out neatly, as summarized in the following theorems. 

Multivector directional derivatives exist.

\begin{gbox}
\begin{thm}[Existence of MDDs]
\label{thm:thm:mdd-exists}
For any metric-compatible affine connection $\Dres$, there exists a multivector directional derivative $D$ that equals $\Dres$ when restricted to act only on vectors.
\end{thm}
\begin{pf}
Proved by construction in Appendix~\ref{sec:app:mdd-existence-proof}.
\end{pf}
\end{gbox}

And they are in bijective correspondence with metric-compatible affine connections.

\begin{gbox}
\begin{thm}[MDDs $\longleftrightarrow$ metric-compatible connections]
\label{thm:thm:mdd-affine-bijection}
The restriction to act on vector fields is a bijection from multivector directional derivatives to metric-compatible affine connections.
\end{thm}
\begin{pf}
The axioms of Definition~\ref{thm:def:directional-deriv} with Theorem~\ref{thm:thm:dd-props-1} ensure that if any $D$ exists the restriction $\Dbar$ of $D$ is a metric-compatible affine connection. Expanding in a canonical multivector frame shows that any $DA$ can be evaluated in terms of $\Dbar$, thus $\Dbar = \Dbar '$ implies $D = D'$ so the map is one-to-one. By Theorem~\ref{thm:thm:mdd-exists}, for any metric-compatible affine connection $\Dres$, there exists an MDD $D$ that restricts to $\Dres$. Thus MDDs exist and the map is onto, so the map is a bijection.
\end{pf}
\end{gbox}
 
Because of this correspondence, multivector directional derivatives can be uniquely specified by a set of metric-compatible connection coefficients. These can be expressed in terms of an arbitrary vector basis $\e_i$ such that
\begin{equation*}
    g_{ij} = \e_i \cd \e_j \, ,
    \qquad \qquad \qquad
    [\e_i, \e_j] = L_{ijk} \, \e^k \, ,
\end{equation*}
with reciprocal basis $\e^i$ such that $\e^i \cd \e_j = \dd^i_j$, as usual.

\begin{gbox}
\begin{cor}[Connection coefficients]
\label{thm:cor:conn-coeffs}
A multivector directional derivative $D$ is uniquely specified by its connection coefficients $\Gamma_{ijk}$ in any vector basis $e_i$, defined by
\begin{equation*}
D_{\e_i} \e_j = \Gamma_{ijk} \, \e^k ,
\end{equation*}
which can be arbitrary other than the restriction
\begin{equation*}
\Gamma_{ijk}+\Gamma_{ikj}=\pp_{\e_i} \, g_{jk} .
\end{equation*}
\end{cor}
\begin{pf}
By Theorem \ref{thm:thm:mdd-affine-bijection}, each metric-compatible affine connection extends uniquely to a multivector directional derivative.
\end{pf}
\end{gbox}

Equivalently, the connection coefficients can be parameterized in terms of the contorsion coefficients, separating out the standard term.

\begin{gbox}
\begin{cor}[Contorsion coefficients]
\label{thm:cor:contorsion-coeffs}
A multivector directional derivative $D$ is uniquely specified by its contorsion coefficients $\chi_{ijk}$ in any vector basis $e_i$, which can be arbitrary other than the restriction
\begin{equation*}
\chi_{ijk}+\chi_{ikj}=0 \, .
\end{equation*}
With $\G_{ijk}$ as in Corollary \ref{thm:cor:conn-coeffs}, the contorsion coefficients $\chi_{ijk}$ are defined by
\begin{equation*}
    \Gamma_{ijk}  
        = \tfrac{1}{2} \left( \pp_{\e_i}g_{jk} - \pp_{\e_k}g_{ij} + \pp_{\e_j}g_{ki} \right)
        + \tfrac{1}{2} \left( L_{ijk} - L_{jki} + L_{kij} \right)
        + \chi_{ijk} \; .
\end{equation*}
Additionally, $D$ is torsion-free if and only if $\chi_{ijk} \equiv 0$.
\end{cor}
\begin{pf}
Torsion for multivector derivatives is defined in the following section. The rest follows immediately from Corollaries \ref{thm:cor:conn-coeffs} and \ref{thm:cor:contorsion}.
\end{pf}
\end{gbox}

It is also sometimes useful to know the directional derivatives of the reciprocal basis. With the following theorem it becomes easy to evaluate derivatives with any combination of the basis and reciprocal basis elements. 

\begin{gbox}
\begin{prop}[Reciprocal connection formula]
\label{thm:prop:reciprocal-coeffs}
Let $D$ be an MDD with connection coefficients defined by $D_{\e_i} \e_j = \G_{ijk} \, \e^k$. Then
\begin{equation*}
    D_{\e_i} \e^j = - (\G_{ilm} \, g^{mj}) \, \e^l
\end{equation*}
gives a formula for the derivatives of the reciprocal basis.
\end{prop}
\begin{pf}
$
D_{\e_i} \e^j = 
D_{\e_i} (g^{jm}\e_m) = 
\pp_{\e_i} (g^{jm}) \, \e_m + g^{jm} D_{\e_i} (\e_m) = 
(\pp_{\e_i} (g^{jm}) g_{ml} + g^{jm} \, \G_{iml}) \e^l = 
(-\pp_{\e_i} (g_{ml}) + \G_{iml}) g^{jm} \, \e^l = 
 - \G_{ilm} g^{jm} \, \e^l
$.
The preceding steps made use of metric compatibility in the form $\pp_{\e_i} g_{jk} = \G_{ijk}+\G_{ikj}$ and used a  product rule expansion of the form $\pp_{\e_i} (g^{jm} g_{mk}) = 0$.
\end{pf}
\end{gbox}

This section has established the existence of a set of natural directional derivative operators on the space of multivector fields. Compared to affine connections, multivector directional derivative operators naturally act on a more useful variety of objects, obey a more intuitive product rule, and have useful properties arising from a minimal set of assumptions.

\section{Torsion-free multivector directional derivative}

We have seen that each multivector directional derivative (MDD) uniquely corresponds to a metric-compatible affine connection. Also like affine connections, MDDs can be characterized by their torsion, and there is a unique torsion-free MDD corresponding to the Levi-Civita affine connection.

The definition of torsion for an MDD is the same as for affine connections (see Definition~\ref{thm:def:metric-compatibility-torsion}). It was already shown that there is a unique metric-compatible torsion-free affine connection. This leads immediately to the following statement.

\begin{gbox}
\begin{prop}[Torsion-free MDD]
\label{thm:prop:del}
There is a unique torsion-free multivector directional derivative with contorsion coefficients
\begin{equation*}
    \chi_{ijk} = 0 
\end{equation*}
which is given the special notation $\DD$ and called the \textit{torsion-free (or Levi-Civita) derivative}.
\end{prop}
\begin{pf}
Existence is guaranteed by Theorem~\ref{thm:thm:mdd-affine-bijection}. Uniqueness and $\chi_{ijk}=0$ follow from Corollary~\ref{thm:cor:contorsion}.
\end{pf}
\end{gbox}

It is convenient to also have a special notation for the torsion-free connection coefficients.

\begin{gbox}
\begin{cor}[Torsion-free connection coefficients]
\label{thm:def:del}
The connection coefficients $\Gbar_{ijk}$ for the torsion-free derivative are defined by $\DD_{\e_i} \e_j = \Gbar_{ijk} \, \e^k$ in an arbitrary basis $\e_i$, and given by 
\begin{equation*}
    \Gbar_{ijk}  
        = \tfrac{1}{2} \left( \pp_{\e_i}g_{jk} - \pp_{\e_k}g_{ij} + \pp_{\e_j}g_{ki} \right)
        + \tfrac{1}{2} \left( L_{ijk} - L_{jki} + L_{kij} \right) \,.
\end{equation*}
These are called the \textit{standard (or Levi-Civita) connection coefficients}.
\end{cor}
\begin{pf}
Set $\chi_{ijk}=0$ in general form of $\G_{ijk}$ (see Corollary~\ref{thm:cor:contorsion-coeffs}).
\end{pf}
\end{gbox}

In standard Riemannian geometry, the induced affine connection on an embedded submanifold of Euclidean $\mathbb{R}^n$ is always metric-compatible and torsion-free \cite{lee97}. Also, geodesics are autoparallels of the metric-compatible torsion-free connection. These facts motivate the usual acceptance of the Levi-Civita connection as the natural choice of connection on arbitrary manifolds. In the present context, all MDDs are metric-compatible, and we identify the unique torsion-free derivative as a natural choice of MDD. The theory of embedded submanifolds has not been explicitly written down yet in the current context, but presumably the same special properties of the torsion-free derivative continue to hold.

Every MDD $D$ can be expressed as the sum of the torsion-free derivative $\DD$ and a contorsion operator $Q$.

\begin{gbox}
\begin{defn}[Contorsion operator]
\label{thm:def:Q}
Let $D$ be an MDD. Then $D$ can be expressed by
\begin{equation*}
    D_a A = \DD_a A + Q_a A
\end{equation*}
Where $Q$ is called the \textit{contorsion operator for $D$}.
\end{defn}
\end{gbox}

The contorsion operator has some useful properties.

\begin{gbox}
\begin{prop}[Contorsion operator properties]
\label{thm:prop:Q-props}
Let $Q = D - \DD$ be the contorsion operator for an MDD $D$. Then
\\[6pt]
\begin{tabular}{r p{.9\textwidth}}
(i)
&
$Q$ has all the properties of Definition~\ref{thm:def:directional-deriv} except that \ref{thm:def:directional-deriv}(ii) is replaced by $Q_a \ph = 0$ for scalar fields. Also $Q$ is grade-preserving.
\\[6pt]
(ii)
&
$Q_{\e_i} \e_j = \chi_{ijk} \, \e^k$.
\\[8pt]
(iii)
&
$Q$ is a tensor field (see Section~\ref{sec:tensor-fields}).
\end{tabular}
\end{prop}
\begin{pf}
(i) All properties follow from direct calculation applying the properties of $D$ and $\DD$.
(ii) Direct calculation.
(iii) $Q$ is a tensor field if it is pointwise linear in both arguments. It is automatically pointwise linear in the direction argument by definition. In the field argument it is pointwise linear since it obeys the product rule and annihilates scalar fields: $Q_a(\ph A)=Q_a(\ph) A + \ph \, Q_a(A) = \ph \, Q_a(A)$.
\end{pf}
\end{gbox}

Since the contorsion operator is a tensor field, it can also be called the \textit{contorsion tensor}. In contrast, neither $D$ nor $\DD$ are tensor fields, since they violate pointwise linearity. Thus, under a change of basis, the $\chi_{ijk}$ transform like tensor coefficients while the $\Gbar_{ijk}$ do not.

In a holonomic orthonormal basis the torsion-free derivative has $\Gamma_{ijk} \equiv 0$. It follows that a manifold is flat (in the usual sense of having zero Riemann curvature and being locally isometric to a Euclidean or Minkowski space of some signature) if and only if there exists a smooth holonomic orthonormal basis.

\section{Tensors}
\label{sec:tensors}

This section generalizes the notion of tensor, from a tensor on a vector space to a tensor on a geometric algebra. The next section will then consider tensor fields on a manifold.

Typically a tensor is defined as a multilinear map from several copies of a vector space (and/or its dual) to $\mathbb{R}$. When a vector space is extended to a geometric algebra, how should the notion of tensor be extended?

There are several ways one could do so. A straightforward approach would be to simply apply the usual definition using the vector subspace of multivectors. This seems not to make enough use of the power of geometric algebra. On the other extreme, one could take a tensor to be a multilinear map of several copies of the geometric algebra back to itself; this would be a bit too general for easy analysis, however, and would undermine the concepts of tensor rank and signature. Each of these could be said to capture the basic idea that a tensor is a multilinear function of multivectors. 

We will use a definition, in the same spirit, which is a middle ground of those extremes; the inputs and output of the multilinear map are each restricted to a fixed-grade subspace of the geometric algebra. Since the geometric algebra approach has identified vectors with their duals, there will be no need to include the dual space.

\begin{gbox}
\begin{defn}[Tensor]
Let $\GA$ be a geometric algebra. Denote by $\grade{\GA}_k$ the linear space of grade-$k$ multivectors in $\GA$. Then a map
$$ T: \grade{\GA}_{k_1} \times \ldots \times \grade{\GA}_{k_N} \to \grade{\GA}_{k_0} $$
is called a \textit{tensor} on $\GA$ if it is \textit{multilinear}, meaning that for every input slot
$$ T(\ldots, \; \alpha A + \beta B, \; \ldots) =  \alpha \, T(\ldots, \; A, \; \ldots) + \beta \, T(\ldots, \; B, \; \ldots)$$
for all scalars $\alpha, \beta$ and valid (correct grade for input slot) multivectors $A,B$. 

If $T$ is a tensor then its \textit{signature} is $(k_1, \ldots, k_N: k_0)$, its \textit{rank} is $k_1 + \ldots + k_N + k_0$, and its \textit{number of inputs} is $N$.
\end{defn}
\end{gbox}

Traditional tensors have, in this formalism, the signature $(1,\ldots,1:0)$. One often finds, however, that the practice of setting $k_0=0$ in the traditional definition can require annoying maneuvering to make multilinear maps into tensors (for example when converting the Riemann curvature function to the Riemann tensor). Allowing more general tensor outputs is one benefit of the multivector formalism.

Tensors can be manipulated in various ways. (Note: to simplify notation in the remainder of the section, assume input multivectors in tensor arguments always have the correct grade.)

\begin{gbox}
\begin{prop}[Tensor operations]
\label{thm:prop:tensor-ops}
Tensors on a geometric algebra $\GA$ admit the following operations:
\\[4pt]
    \begin{tabular}{r p{.9\textwidth}}
    (i)
    &
    \textit{Addition.}
    $$ (T+S)(A,\ldots,B) = T(A,\ldots,B) + S(A,\ldots,B) \, . $$
    \\[-12pt]    
    (ii)
    &
    \textit{Scalar multiplication.}
    $$ (\alpha \, T)(A,\ldots,B) = \alpha \, T(A,\ldots,B). $$
    \\[-12pt]
    (iii)
    &
    \textit{Multiplication.} For tensors with scalar outputs (signatures $(\ldots:0)$):
    $$(T \otimes S)(A,\ldots,B,C,\ldots,D) = T(A,\ldots,B) \,S(C,\ldots,D) \, .$$
    This is called the \textit{tensor product}.
    \\[6pt]
    (iv)
    &
    \textit{Contraction.} A tensor of signature $(\ldots,1,\ldots,1,\ldots,:k_0)$ (that is, a tensor with at least two grade-$1$ inputs) can be contracted:
    $$ \bar{T}(\ldots,\ldots,\ldots) = \sum_i \; T(\ldots,\e_i,\ldots,\e^i,\ldots) \, , $$
    where $\e_i$ is an arbitrary basis. The contraction is well-defined, since it can be shown to be independent of the evaluation basis. The result is a tensor of signature $(\ldots,\ldots,\ldots,:k_0)$ (two fewer inputs and rank reduced by $2$).
    \\[0pt]
    \end{tabular}
\end{prop}
\begin{pf}
(i-iii) Trivial. (iv) $\bar{T}$ is independent of the evaluation basis. Let $\e_i$ and $\E_i$ be any two bases. Suppressing unnecessary arguments, we now show that $\bar{T}=T(\e_i,\e^i)=T(\E_i,\E^i)$ (with summation convention). The proof is direct: $T(\e_i,\e^i) = T((\e_i \cd \E^j)\E_j,(\e^i \cd \E_k)\E^k) = (\e_i \cd \E^j) (\e^i \cd \E_k) \; T(\E_j,\E^k) = (\E^j \cd \E_k) \; T(\E_j,\E^k)  = {\dd^j}_k \; T(\E_j,\E^k) = T(\E_j,\E^j)$. The proof also confirms that the order of the upper and lower index are unimportant.
\end{pf}
\end{gbox}

The zero function $Z(A,\ldots,B)=0$ is a valid tensor (the \textit{zero tensor}) for any signature ($k_0$~can be arbitrary since the zero multivector counts as every grade). Given the zero tensor, it can be shown that the space of tensors of a fixed signature forms a finite-dimensional linear space under addition and scalar multiplication. It is implicit in the above notation that the operations of addition and scalar multiplication each operate within a space of fixed signature, while multiplication and contraction operate in the space of all tensors of arbitrary signature.

It is sometimes useful to express tensors in terms of components. In general, the components of an arbitrary tensor relative to the multivector basis $\e_J$ (with reciprocal basis $\e^J$) are defined by
\begin{equation}
    T(\e_{J_1}, \ldots, \e_{J_N}) = T_{{J_1} \, \ldots \, {J_N} {J_0}} \;\, \e^{J_0} \, ,
\end{equation}
where the multivector basis elements are evaluated only at the grade appropriate for each slot. Components for the same tensor can also be evaluated in the reciprocal basis,
\begin{equation}
    T(\e^{J_1}, \ldots, \e^{J_N}) = T^{{J_1} \, \ldots \, {J_N} {J_0}} \;\, \e_{J_0} \, ,
\end{equation}
which are written with an upper index. Alternately, mixed up-down tensor components can be obtained by choosing to evaluate each slot using the basis or reciprocal basis. Vector basis indices (corresponding to grade-1 slots) can be raised and lowered with the metric, as seen below, but general multivector basis indices cannot.

When higher grade multivectors are involved, the component notation is a little bit clumsy. Each multivector index $J$ is equivalent to a sequence of several vector indices $j$, so that, written out fully in terms of vector indices, the number of vector indices is equal to the rank of $T$. But since not every sequence of vector basis elements is represented in the multivector basis, the rank alone does not uniquely determine the number of independent components of a tensor; this number can, however, be determined from the signature. On the other hand, if all indices are vector basis indices, the total number of independent components is $d^{\,r}$, where $r$ is the rank and $d$ is the dimension of the space of vectors.

Components are most useful when working with traditional tensors (of signature $(1,\ldots,1:0)$). In that case,
\begin{equation}
\label{eqn:tensor-components}
    T(\e_{j_1}, \ldots, \e_{j_N}) = T_{{j_1} \, \ldots \, {j_N} } \, ,
\end{equation}
and likewise for the upper index components in the reciprocal basis, or mixed components using both bases. For a vector basis one has $\e^i = g^{ij} \e_j$ so that, by linearity, component indices can be raised or lowered with the metric, for example as in
\begin{equation}
    \begin{array}{rcl}
        {T^{i}}_j   & = & T(\e^i,\e_j) \;=\; T(g^{ik} \e_k, \e_j) \;=\; g^{ik}\, T( \e_k, \e_j) \\[4pt]
         & = & g^{ik} \, T_{kj}  \, .
    \end{array}
\end{equation}
Similarly, the basis transformation formula for components is derived by linearity, from
\begin{equation}
    \begin{array}{rcl}
        T(\e_{j_1}, \ldots, \e_{j_N})  & = & T((\e_{j_1}\cd \E^{k_1})\E_{k_1}, \ldots, (\e_{j_N}\cd \E^{k_N})\E_{k_N}) \\[4pt]
         & = & (\e_{j_1}\cd \E^{k_1}) \,\ldots \, (\e_{j_N}\cd \E^{k_N}) \; T(\E_{k_1}, \ldots, \E_{k_N}) \, ,
    \end{array}
\end{equation}
for arbitrary vector bases $\e_i$ and $\E_i$. Corresponding formulas in terms of upper or mixed components can be extrapolated straightforwardly. Moreover, basically the same formalism holds for tensors of signature $(1,\ldots,1:1)$ --- as long as only vector basis indices are involved, indices behave the same as in traditional treatments.

Every vector can be naturally identified with a certain tensor, which we call its \textit{tensor conjugate}. For any vector $b$, define the tensor $\Hat{b}$ acting on a vector input $a$ by
\begin{equation}
    \Hat{b} \,(a) = b \cd a \, .
\end{equation}
This $\Hat{b}$ tensor is the object usually called the \textit{one-form} conjugate to the vector $b$. In the case of reciprocal basis frames, their properties ensure that as tensors they obey
\begin{equation}
\label{eqn:basis-oneforms}
    \Hat{\e}^i \,(\e_j) = \e^i \cd \e_j = \dd^i_j \, ,
\end{equation}
a property which will be useful below.

More generally, we can provide a tensor conjugate to any fixed-grade multivector. Let $B_k$ be a multivector of grade $k$. The tensor $\Hat{B_k}$ conjugate to $B_k$ is a tensor mapping $k$ vector inputs to a scalar (signature $(1,\ldots,1:0)$ with $k$ input slots), defined by
\begin{equation}
    \label{eqn:tensor-conjugate}
    \Hat{B_k} \,(a_1,\ldots,a_k) = B_k \cd (a_1 \w \ldots \w a_k) \, .
\end{equation}
This is the usual tensor representation of a $k$-form. There is also an alternate way to map grade-$k$ multivectors to tensors, defined by
\begin{equation}
\label{eqn:tensor-conjugate-2}
    \Breve{B_k} \,(A)= B_k \cd \grade{A}_k \, ,
\end{equation}
resulting in a tensor of signature $(k:0)$. In both (\ref{eqn:tensor-conjugate}) and (\ref{eqn:tensor-conjugate-2}) the grade index would not usually be made explicit, but it is given in these cases to clarify how the grade affects the definitions. Clearly, there is a close relationship between the two possibilities, and it seems likely that each could prove useful in different situations. Although (\ref{eqn:tensor-conjugate-2}) is more in line with the spirit of the multivector approach, we take (\ref{eqn:tensor-conjugate}) as the standard definition in order to make closer contact with standard tensor calculus.

One can multiply the tensor conjugates $\hat{a},\hat{b}$ of vectors $a,b$ via the tensor product, giving
\begin{equation}
    (\Hat{a}\otimes \Hat{b}) \; (c,d) = (a \cd c) (b \cd d),
\end{equation}
with components
\begin{equation}
    (\Hat{a}\otimes \Hat{b}) \; (\e_i,\e_j) = (a \cd \e_i) (b \cd \e_j) = a_i \, b_j \; .
\end{equation}
Thus arbitrary tensors of the traditional type (signature $(1,\ldots,1:0)$) be built up in terms of the tensor product as
\begin{equation}
\label{eqn:general-tensor-product}
    T = T_{j_1 \, \ldots \, j_N} \;\; \Hat{\e}^{\, j_1}\otimes \ldots \otimes \Hat{\e}^{\,j_N} .
\end{equation}
This is a completely general tensor of signature $(1,\ldots,1:0)$. Applying the definition of the tensor product along with (\ref{eqn:basis-oneforms}) shows that components in (\ref{eqn:general-tensor-product}) are consistent with (\ref{eqn:tensor-components}). While this form of tensor notation is sometimes useful, it also tends to obscure the basic nature of a tensor as a function, and here we avoid it when possible.

This section has defined tensors on a geometric algebra, studied some of their basic properties, and connected the present definitions to traditional ones. Next this is extended to study tensor fields on a manifold.

\section{Tensor fields}
\label{sec:tensor-fields}

In traditional differential geometry, and especially in general relativity, tensor fields are given a role of utmost significance. Even vectors, the most basic objects of the theory, are usually treated as a special case of tensors. This leads to a development of tensor calculus in which the role of vectors as geometric objects and the role of tensors as linear maps are mixed together. As a result, the natural product rule for multivectors (see above) and natural chain rule for tensors (see below) each become obscured. The geometric algebra approach allows somewhat more conceptual clarity, by separating the role of multivectors (geometric objects), from the role of tensors (linear maps). This helps facilitate a simple derivation of the tensor covariant derivative from the chain rule.

Recall from the previous section that a tensor is defined as a linear operator on a geometric algebra. In the case of a tensor field on a manifold, the underlying geometric algebra at each point is the geometric tangent space $GT_p M$. Tensor fields are defined only when the tensor at every point has the same signature, as follows. (Note that the definition is rather informal. We avoid a more formal definition, to avoid introducing the language of bundle sections, as it would not contribute much to the discussion.)

\begin{gbox}
\begin{defn}[Tensor field]
A \textit{tensor field} $T$ of signature $\mathcal{S} =(k_1, \ldots, k_N: k_0)$ is an assignment of a tensor of signature $\mathcal{S}$ on $GT_p M$ to each point on (or in a neighborhood on)~$M$. A tensor field is called \textit{smooth} if smooth inputs always lead to a smooth output.
\end{defn}
\end{gbox}

Each point is assigned its own multilinear function. It follows that the multilinearity property of a tensor field is evaluated pointwise, so that even for a non-constant scalar field $\ph$, the tensor field $T$ has the pointwise linearity $T(\ph A) = \ph \, T(A)$. A map like $T$ which has linearity with respect to constant scalars, but not to non-constant scalar fields, is not a tensor field (the multivector directional derivative $D_a$ in the direction of a fixed vector field $a$, for example, is not a tensor field for precisely this reason).

Tensor fields admit the same operations as tensors (see Proposition \ref{thm:prop:tensor-ops}), with each operation being evaluated pointwise. Tensor fields also admit the additional operation of local scalar multiplication (multiplication by a non-constant scalar field); the difference between this and global scalar multiplication should be clear from context. The space of tensor fields of a fixed signature forms an infinite-dimensional linear space under addition and global scalar multiplication. The basic definitions and operations for tensor fields have now been established.

How is the directional derivative of a tensor field to be taken? Since the multivector directional derivative allows differentiation of the inputs and output of a tensor field, intuition suggests that some kind of chain rule should apply. Let us investigate this further by first considering a linear function of some real variables.

Consider a function $f(x,g(x),h(x))$. This function depends directly on the real number $x$, as well as on two quantities $g$ and $h$ that each vary with $x$. As $x$ varies, the total differential is
\begin{equation}
    df = \tfrac{\pp f}{\pp x} \, dx + \tfrac{\pp f}{\pp g} \, dg + \tfrac{\pp f}{\pp h} \, dh \; .
\end{equation}
Suppose now that $f$ is linear in the $g$ and $h$ arguments. Then 
\begin{equation}
    \begin{array}{rcl}
         \tfrac{\pp f}{\pp g} \, dg &=& f(x, g+dg, h) - f(x, g, h) \\[4pt]
                                    &=& f(x,dg, h)
    \end{array}
\end{equation}
by linearity, and likewise for the $h$ term, so that
\begin{equation}
    df = \tfrac{\pp f}{\pp x} \, dx + f(x,dg, h) + f(x,g, dh) \; .
\end{equation}
The lefthand side is the change in the \textit{output} of $f$, while the two rightmost terms represent the change in output due to changing \textit{inputs}. Only the $\frac{\pp f}{\pp x}$ term represents the change in ``$f$ itself"---this is the term that corresponds to the tensor derivative $DT$ below.

Analogously, a tensor field depends both directly on the point of evaluation, and linearly on some multivector fields inputs. A correspondence can be drawn as follows
\begin{equation*}
    \begin{array}{ccccclc}
        df &= &\tfrac{\pp f}{\pp x} \, dx &+ &f(x,dg, h) &+ &f(x,g, dh) \\[4pt]
        \downarrow & &\downarrow & &\downarrow & &\downarrow \\[4pt]
        D_a(T(A,\ldots,B)) &= &(DT)(a,A,\ldots,B)  &+ &T(D_a A,\ldots,B) &+\ldots + &T(A,\ldots,D_a B) .
    \end{array}
\end{equation*}
Note that every term except the $DT$ term is already defined by virtue of the multivector directional derivative. Therefore the chain rule can be used to \textit{define} the quantity $DT$. After giving the formal definition below, we will confirm that such a $DT$ is in fact a well-defined tensor field, and that it is equivalent to the usual tensor covariant derivative.

The tensor derivative of a tensor field $T$ is a new tensor field $DT$ with one additional vector input (corresponding to the derivative's direction argument).

\begin{gbox}
\begin{defn}[Tensor derivative]
Let $T$ be a tensor field of signature $(k_1, \ldots, k_N: k_0)$, and let $D$ be a multivector directional derivative. Define $DT$ by
\begin{equation}
    (DT)(a,A,\ldots,B)  = D_a(T(A,\ldots,B)) -  T(D_a A,\ldots,B) -\ldots - T(A,\ldots,D_a B) \, .
\end{equation}
$DT$ is a tensor field of signature $(1,k_1, \ldots, k_N: k_0)$, called the \textit{tensor derivative} of $T$.
\end{defn}
\begin{pf}
For $DT$ to be a tensor field, it must be pointwise linear in every argument. This can be shown by direct calculation using the pointwise linearity of $T$ and the multivector directional derivative properties. The key step is to note that $D_a(T(\alpha A)) - T(D_a(\alpha A)) = \alpha \, D_a(T(A)) - \alpha \, T(D_a A)$ by cancellation of the unwanted $(\pp_a \alpha)\, T(A)$ terms. 
\end{pf}
\end{gbox}

It is not trivial that, in addition to obeying the correct chain rule, $DT$ is a well-defined tensor field. Some other seemingly reasonable possible definitions do not have this property. For example if one assumed the incorrect definition $(DT)(a,A,\ldots,B)= D_a (T(A,\ldots,B))$, then $DT$ would fail to be a tensor field by violating pointwise linearity.

Sometimes it is useful to use a tensor derivative notation more similar to the usual directional derivative notation.

\begin{gbox}
\begin{defn}[Tensor derivative notation]
\label{thm:def:tensor-derivative-notation}
Let $T$ be a tensor field. As an alternative notation for the tensor derivative, define
\begin{equation}
    \begin{array}{rcl}
        D_a T(A,\ldots,B) &\equiv& (DT)(a,A,\ldots,B)\\[4pt]
                          &\neq  & D_a ( T(A,\ldots,B))
    \end{array}
\end{equation}
In other words, if $T$ is a tensor field, then $D_a T$ generically refers to the tensor derivative of $T$, as opposed to the multivector directional derivative of the output of $T$. 

This notation is convenient but has a definite risk of confusion, so to specify the derivative of the output explicit parentheses should be used and the meaning should be made clear in context.
\end{defn}
\end{gbox}

Likewise, $\DD T$ or $\DD_a T$ denotes the tensor derivative evaluated using $\DD$ (the torsion-free MDD).

This tensor derivative is equivalent to the usual covariant derivative of tensors. To see why, note that the tensor derivative acts in the usual way on components. Consider for example a tensor field
\begin{equation}
    T = T_{ij} \;\; \Hat{\e}^i \otimes \Hat{\e}^j
\end{equation}
and its derivative
\begin{equation}
    DT = D_i T_{jk}  \;\;  \Hat{\e}^i \otimes \Hat{\e}^j \otimes \Hat{\e}^k \, .
\end{equation}
A direct calculation from the definiton of $DT$ reveals
\begin{equation}
\label{eqn:DT-lower-indices}
    D_i T_{jk} = \pp_i T_{jk} - (\G_{ijm} \, g^{ml}) \, T_{lk}  - (\G_{ikm} \, g^{ml}) \, T_{jl}
\end{equation}
which can quickly be checked to be equivalent to the usual expression. Since $DT$ is a tensor, it also follows that
\begin{equation}
\label{eqn:DT-upper-indices}
\begin{array}{rcl}
    D_i T^{jk}  &=& g^{mj} g^{nk} \, D_i T_{mn} 
    \\[8pt]
                &=& \pp_i T^{jk} + (\G_{ilm} \, g^{mj}) \, T^{lk}  + (\G_{ilm} \, g^{mk}) \, T^{jl}
\end{array}
\end{equation}
where the second step uses metric-compatibility. The same expression for $D_i T^{jk}$ can also be confirmed directly from the definition of $DT$ using an application of Proposition~\ref{thm:prop:reciprocal-coeffs}. Note that $D_i T^{jk}$ is merely a shorthand for the mixed-index components ${(DT)_i}^{jk}$ of the (invariantly defined) tensor derivative $DT$. The shorthand helps make clear which index corresponds to the derivative direction. The lower-index expression is derived more straightforwardly because $\G_{ijk}$ was defined in terms of lower indices.

The expressions (\ref{eqn:DT-lower-indices}--\ref{eqn:DT-upper-indices}) for the components of the tensor derivative can be naively generalized to tensors of arbitrary rank and mixed upper/lower indices. Note the equivalence of both expressions to the usual ones where $\G$ is defined with one upper and two lower indices. The fact that the tensor derivative as defined here gives the usual expression in components is sufficient to prove its equivalence to the usual covariant derivative of tensors. Alternatively, one can check that it has all the defining properties usually ascribed to the covariant derivative.

\begin{gbox}
\begin{thm}
\label{thm:thm:dt-equiv-covariant}
$DT$ is equivalent to the usual covariant derivative of tensor fields (see e.g.~\cite{lee97}).
\end{thm}
\begin{pf}
This can be proved in two ways, either axiomatically, or in components. It was already shown in Equations~(\ref{eqn:DT-lower-indices}--\ref{eqn:DT-upper-indices}) and surrounding text that the components of $DT$ are equal to the components of the usual covariant deriative. (Note that in a holonomic basis our $\Gbar_{ijm} g^{mk}$ reduce to the usual Levi-Civita connection coefficients.) This is sufficient to prove equivalence. On the other hand it can also be shown that $DT$ acts as usual on scalar and vector fields, obeys the product rule (\ref{eqn:dt-product-rule}), and commutes with contractions~(\ref{eqn:dt-commutes-contractions}), forming an axiomatic proof. The proof of the product rule is straightforward from the definition. To prove the contraction formula note that Proposition~\ref{thm:prop:reciprocal-coeffs} implies \mbox{$T(D_{\e_i} \e^k, \e_k) + T( \e^k, D_{\e_i} \e_k) = 0$} so the extraneous term in $\overline{D_a T}$ vanishes. Note that metric compatibility is essential for the proof that contraction commutes.
\end{pf}
\end{gbox}

On tensors which correspond directly to multivectors, the tensor derivative and multivector derivative are equivalent. In particular, using dot and wedge compatibility of the multivector derivative, it is straightforward to show that for any tensor conjugate $\hat{A}$ of a grade-$k$ multivector~$A$, the tensor and multivector derivatives are related by
\begin{equation}
    D_a \, \hat{A} = \widehat{D_a A}
\end{equation}
with each side of the expression being a tensor of rank $k+1$ if $a$ is regarded as a free input variable.

Additional useful properties are that the tensor derivative obeys the product rule
\begin{equation}
\label{eqn:dt-product-rule}
    D_a (T \otimes S) = D_a T \otimes S + T \otimes D_a S
\end{equation}
over tensor products and commutes with contractions in the sense that with the contraction $\bar{T}$ of~$T$ defined as in Proposition~\ref{thm:prop:tensor-ops}
\begin{equation}
\label{eqn:dt-commutes-contractions}
    \overline{D_a T} = D_a \bar{T}
\end{equation}
assuming the same arguments are contracted on each side. (See Theorem~\ref{thm:thm:dt-equiv-covariant} for proofs.)

There is no particularly good reason to think of the metric as a tensor, rather than just as a set of coefficients determining the dot product. Nonetheless, if one insists, they may define the metric tensor by $\hat{g}(a,b)=a \cd b$. It then follows from metric compatibility of the multivector directional derivative that the tensor derivative $D\hat{g} = 0$ vanishes. This is the sense in which the tensor derivative is compatible with the metric tensor.

\section{Gradient, divergence, curl}

Geometric calculus admits intuitive generalizations of the gradient, divergence, and curl of vector calculus~\cite{macdonald2017}. Here we generalize these notions to the present formalism.

The gradient operator is defined as follows.

\begin{gbox}
\begin{defn}[Gradient]
Let $D$ be a multivector directional derivative, and $A$ be a multivector field. The \textit{gradient} operator $\D$ is defined by
$$ \D A = \e^i \, D_{\e_i} A \; ,$$
where $\e_i$ is an arbitrary basis. $\D A$ is well-defined, since it can be shown that $\e^i \, D_{\e_i} A = \E^j \, D_{\E_j} A$ if $\E_i$ is any other basis. There is an implicit sum in this definition due to the repeated index.
\end{defn}
\begin{pf}
Expanding $\E_i$ in the $\e_i$ basis gives
$
\E^i \, D_{\E_i} A
=
(\E^i \cd \e_j) \, \e^j \; D_{(\E_i \cd \e^k) \, \e_k} A
=
(\E^i \cd \e_j) \, (\E_i \cd \e^k) \, \e^j \; D_{ \e_k} A
=
(\e^k \cd \e_j) \, \e^j \; D_{ \e_k} A
=
\e^k \, D_{\e_k} A
$.
Thus $\D$ is well-defined and independent of basis.
\end{pf}
\end{gbox}

It is useful to think of the gradient operator as a vector $\D = \e^i \, D_{\e_i}$, where all vector aspects are included in $\e^i$, and $D_{\e_i}$ is a scalar (grade-preserving) operator acting to the right. For this reason $\D$ is often called a \textit{vector operator}. Note that $\e^i$ and $D_{\e_i}$ do not commute, so the order of the two terms is significant. All properties of $\D$ are derived from the multivector directional derivative properties explored earlier. Certain properties are explored further in Appendix~\ref{sec:interp-grad}.

Since by (\ref{eqn:ext-fund-iden}) vectors obey $aB = a \cd B + a \w B$, the gradient can be decomposed into a grade-lowering and grade-raising term corresponding to the divergence and curl.

\begin{gbox}
\begin{thm}[Gradient $=$ divergence $+$ curl]
Let $D$ be a multivector directional derivative, let $\e_i$ be an arbitrary basis, and let $A$ be a multivector field. Define the \textit{divergence} and \textit{curl} by
\begin{equation*}
    \begin{array}{ll}
        \D \wcd A = \e^i \wcd D_{\e_i} A,  \qquad \qquad & \textrm{(divergence)} \\[4pt]
        \D \w   A = \e^i \w   D_{\e_i} A.  \qquad \qquad & \textrm{(curl)}
    \end{array}
\end{equation*}
Then the gradient equals the divergence plus the curl,
\begin{equation}
    \D A = \D \cd A + \D \w A \,.
\end{equation}
These definitions are consistent with the usual definitions of dot and wedge product, and can be shown to be basis-independent and thus well-defined.
\end{thm}
\begin{pf}
The proof of basis-independence is identical to the proof for the gradient, as it makes use only of linearity. The main theorem follows from (\ref{eqn:ext-fund-iden}) since $\D A = \e^i D_{\e_i} A = \e^i \cd D_{\e_i} A + \e^i \w D_{\e_i} A =  \D \cd A + \D \w A$.
\end{pf}
\end{gbox}

Some other related notations will also prove useful.

\begin{gbox}
\begin{defn}
To avoid ambiguity, here are some notational conventions:
\begin{equation*}
\begin{array}{rcl}
    \D^{\,2}  A &=& \D(\D A), \\[4pt]
    (a \cd \D) \, A &=& D_{a} A, \\[4pt]
    (\D^{\,2})\,  A &=& (\e^i \, \e^j) \; D_{\e_i}  D_{\e_j}  A, \\[4pt]
    (\D \wcd \D) \, A &=& (\e^i \wcd \e^j) \; D_{\e_i}  D_{\e_j}  A, \\[4pt]
    (\D \w  \D) \, A &=& (\e^i \w \e^j) \; D_{\e_i}  D_{\e_j}  A \\[4pt]
    (\D \w  \D) \cd A &=& (\e^i \w \e^j) \cd (D_{\e_i}  D_{\e_j}  A) 
\end{array}
\end{equation*}
Similar definitions for $\D^n A$ and $(\D^n)A$ and $(\D \w \ldots \w \D) A$, as well as for expressions like $(\D \cd \D) \w A$ and its cousins, can be extrapolated straightforwardly from those above. Note the importance of explicit parenthesis in defining the notation.

These notations can be checked to be consistent with the definitions of dot and wedge, and the ones depending on a basis can be checked to be well-defined.
\end{defn}
\end{gbox}

As noted earlier, the torsion-free directional derivative is denoted by $\DD$. Therefore take
\begin{equation*}
    \DDD A = \e^i \, \DD_{\e_i} A
\end{equation*}
to express the gradient associated with the torsion-free derivative, and likewise for its associated divergence, curl, and the rest of the notation. The next section shows that the torsion-free curl associated with this gradient is closely related to the exterior derivative of differential forms. The relation of $\DDD \w$ to the usual curl from vector calculus, as well as some more properties of the gradient, are also given at the end of the next section.

\section{Differential forms, exterior derivative}
\label{sec:forms}

The space of differential forms (totally antisymmetric tensors) has an equivalent structure to the exterior subalgebra of multivector fields. In particular, the tensor conjugate map $A \mapsto \hat{A}$ defined in Section~\ref{sec:tensors} provides a bijection from multivectors to forms. Here we investigate this correspondence in more detail, and show that the theory of differential forms is fully included within the present methods.

For reasons to be justified shortly, define the exterior derivative $d = \DDD \w$ as the torsion-free curl.

\begin{gbox}
\begin{defn}[Exterior derivative]
Define the \textit{exterior derivative} $d$ of a multivector $A$ by
\begin{equation*}
    d A = \DDD \w A \; .
\end{equation*}
That is, the exterior derivative $d$ is equivalent to the torsion-free curl.
\end{defn}
\end{gbox}

When acting on scalar fields $\ph$ the exterior derivative
\begin{equation}
    d \ph = \DDD \w \ph = \DDD \ph = e^{i} \, \pp_{\e_i} \ph
\end{equation}
is merely the usual gradient. More generally, the exterior derivative has a simple expression when general multivectors are expressed in a gradient basis. In particular, in terms of the coordinate gradient basis $dx^i$ for a coordinate system $x^i$ the exterior derivative of an arbitrary multivector field $A = \sum_J A_J \, dx^J$ evaluates to
\begin{equation}
    dA = \sum_J dA_J \w dx^J \, 
\end{equation}
where as usual $dx^J$ is the multivector basis constructed by wedging the vector basis $dx^i$. This formula is proved within the proof of Corollary~\ref{thm:cor:d-props-2} below. In any basis which is not a gradient basis the expression of $dA$ does not have such a simple form.

The exterior derivative has some important properties.

\begin{gbox}
\begin{thm}[Exterior derivative properties]
\label{thm:thm:d-props}
The exterior derivative $d=\DDD \w$ of multivectors has the following properties:\\[4pt]
\begin{tabular}{r p{.9\textwidth}}
(i) & $d(A + B) = dA + dB$. 
\\[6pt]
(ii) & If $\ph$ is a scalar field then $a \cd d \ph = \pp_a \ph$ for all vector fields $a$. 
\\[6pt]
(iii) & If $\ph$ is a scalar field then $d^2 \ph \equiv d(d \ph)=0$. 
\\[6pt]
(iv) & If $A_j$ and $B_k$ are multivectors of fixed grades $j$ and $k$ respectively, then
\begin{equation*}
    d(A_j \w B_k) = d(A_j) \w B_k + (-1)^{j} \, A_j \w d(B_k) \, .
\end{equation*}
\\[-20pt]
\end{tabular}
\end{thm}
\begin{pf}
(i) Trivial. 
(ii) $a \cd d\ph = a \cd \e^i \, \pp_{\e_i} \ph = \pp_a \ph$.
(iii) Choose to work in a holonomic coordinate basis $\e_i = e(x_i)$ so that $L_{ijk}=0$. Note that in this basis $e^i \w \DD_{\e_i} e^j = -\Gbar_{ilm} g^{mj} (\e^i \w \e^l) = 0$ since $(\e^i \w \e^l)=-(\e^l \w \e^i)$ and by torsion-freeness $\Gbar_{ilm}-\Gbar_{lim}=L_{ilm}=0$. Similarly note that in this basis $(\e^i \w \e^j) \pp_{\e_i} \pp_{\e_j} \ph = 0$ since coordinate partial derivatives commute. Thus
$d^2 \ph = \e^i \w \DD_{\e_i} ( e^j \pp_{\e_j} \ph)
 = (\e^i \w \DD_{\e_i} e^j) \pp_{\e_j} \ph + (\e^i \w e^j )  \pp_{\e_i} \pp_{\e_j} \ph) = 0$.
 (iv) Note that $\e^i \w A_j = (-1)^j (A_j \w e^i)$ since the wedge product commutes under vector swaps. Thus (suppressing grade labels)
 $d(A \w B) = \e^i \w ( \DD_{\e_i} A \w B + A \w \DD_{\e_i}B)
  = d A \w B + \e^i \w A \w \DD_{\e_i}B 
  = d A \w B + (-1)^j A \w dB $.
\end{pf}
\end{gbox}

Note that the torsion-free assumption was essential in the proof of (iii). Also, together these conditions are sufficient to prove an additional property.

\begin{gbox}
\begin{cor}[Exterior derivative properties]
\label{thm:cor:d-props-2}
The above properties also imply that \\[6pt]
\begin{tabular}{r p{.9\textwidth}}
(v) &  On any multivector field $A$
\begin{equation*}
    d^2 A = \DDD \w (\DDD \w A) = 0 \; .
\end{equation*}
\\[-20pt]
\end{tabular}
\end{cor}
\begin{pf}
Choose to work in the gradient basis $dx^i \equiv \DDD \w x^i$ associated with a coordinate system $x^i$ (see Section~\ref{sec:GM}). Any $A$ can be expanded $A = \sum_J A_J \w dx^J$ where $A_J$ are scalar components. Then 
$dA = \sum_J (dA_J \w dx^J + A_J \w d(dx^J) )$ by properties (i,iv). But $d(dx^J) = d(dx^{j_1} \w \ldots \w dx^{j_N}) = 0$ using properties (iv,iii) since each $dx^i$ has the form $d\ph$. Thus $dA = \sum_J (dA_J \w dx^J)$. Applying the same logic again, one then finds $d^2 A = \sum_J (d^2 A_J \w dx^J)=0$, since $A_J$ are scalars as noted before.
\end{pf}
\end{gbox}

The properties (i-v) are precisely the defining properties of the exterior derivative in differential forms! (See, e.g., \cite{lee02} for a standard treatment.) Conceptually, this is sufficient to make the identification between $d = \DDD \w$ and the usual exterior derivative of differential forms.

This equivalence is made rigorous by the following theorem.
 
\begin{gbox}
\begin{thm}[Equivalence of forms to multivectors]
\label{thm:thm:equivalence-multivectors-forms}
The tensor conjugate map $A \mapsto \hat{A}$ (see Section~\ref{sec:tensors}) is a bijection from multivectors to differential forms which preserves the linear, exterior, and differential structures, in the sense that
\\[8pt]
\begin{tabular}{rl}
    (i) & $\reallywidehat{\alpha A + \beta B} = \alpha \hat{A} + \beta \hat{B}$.\\[4pt]
    (ii) & $\widehat{A \w B} = \hat{A} \w \hat{B}$, \\[4pt]
    (iii) & $\widehat{dA} = \hat{d} \hat{A}$.\\
\end{tabular}
\\[4pt]
Note that $\hat{d}$ is the exterior derivative of differential forms, while $d$ is the exterior derivative of multivectors. Differential form definitions may be found in~\cite{lee02}.
\end{thm}
\begin{pf}
Since forms, unlike tensors, technically allow the formal sum of different signatures, the tensor conjugate mapping for present purposes should be extended to $\hat{A} = \sum_k \widehat{\grade{A}_k}$, with the $k$-form $\widehat{A_k} (a_1, \ldots, a_k) = A_k \cd(a_1 \w \ldots \w a_k)$ as defined previously, where we've replaced $A_k \equiv \grade{A}_k$ for notational convenience. Clearly the image of $A_k$ is a totally antisymmetric tensor, so the image of $A$ is a form. It is trivial to show (i) using linearity of the grade operator and dot product. To proceed, introduce a coordinate system $x^i$ with coordinate basis $\e(x_i)$ and reciprocal basis $dx^i$ in the multivector sense. The image $\wh{dx^i}$ is equivalent to the differential form $\dhat x^i$ since 
$\wh{dx^i}(e(x_j)) = dx^i \cd e(x_j) = \dd^i_j = \dhat x^i(e(x_j))$. 
Now let $dx^{J} = dx^{j_r} \w \ldots \w dx^{j_1}$. Using Equation~(4.12) of \cite{hestenes87} and Proposition 14.11e of \cite{lee02}, it follows that
$
\wh{dx^J}(a_1,\ldots,a_r) 
   = ( dx^{j_r} \w \ldots \w dx^{j_1} ) \cd (a_{k_1} \w \ldots \w a_{k_r})
   = \det(dx^j \cd a_k)
   = \det(\dhat x^j ( a_k))
   = \dhat x^J (a_1,\ldots,a_r).
$
In other words, the multivector basis elements $dx^J$ map to the form basis $\dhat x^J$. Thus, using linearity, if $A = \sum_K A_K dx^K$ then $\hat{A} = \sum_K A_K \dhat x^K$. This formula makes it clear that the map is a bijection, since  $\hat{A}=0$ implies $A_K=0$ implies $A=0$, and the form with arbitrary coefficients $A_K$ is mapped to by the multivector with equal coefficients. It remains to show (ii) and (iii). For (ii), the same logic as above shows that $\reallywidehat{dx^J \w dx^K} = \dhat x^J \w \dhat x^K$. Thus $\widehat{A \w B} = \sum_{JK} A_J B_K \reallywidehat{dx^J \w dx^K} = \hat{A} \w \hat{B}$. Then (iii) also follows, since 
$
\wh{dA} = \sum_K \reallywidehat{dA_K \w dx^K}
    = \sum_K \wh{dA_K} \w \dhat x^K
    = \sum_K \pp_i A_K \, \dhat x^i \w \dhat x^K
    = \dhat \hat{A}
$
using Theorem~14.24 of \cite{lee02} for the final step.
\end{pf}
\end{gbox}

This shows that the theory of multivector fields includes the the theory of differential forms as a part of a more general theory.

A corollary of this equivalence is that the multivector exterior derivative $d = \DDD \w$ is independent of the metric, since it is equivalent to the exterior derivative of differential forms (which is defined without reference to a metric). This fact is somewhat surprising, since the directional derivative $\DD$ underlying $d$ is metric-compatible.

Geometric algebra also has an equivalent of the Hodge dual of differential forms: the \textit{multivector dual}. The dual $A^*$ of a multivector $A$ is defined by $A^* = A I^{-1}$, where $I$ is an oriented unit pseudoscalar~\cite{macdonald2017}. The dual of a $k$-vector $A_k$ is an $(N-k)$-vector representing the orthogonal complement of $A_k$, where $N$ is the total dimension of the space. The multivector dual has the important properties $(A \cd B)^* = A \w B^*$ and $(A \w B)^* = A \cd B^*$ and $0^* = 0$~\cite{macdonald2017}. The dual operation can be used to investigate in more detail the relation between geometric calculus and differential forms. 

The quantity $\DDD \times a$ which in three-dimensional vector calculus is usually called the curl of a vector field $a$ is related to the multivector curl by duality: $\DDD \times a = (\DDD \w a)^*$. In three dimensions this returns a vector orthogonal to the plane of the bivector $\DDD \w a$.

A useful consequence of the properties of duality is that $d^2 = 0$ implies $\DDD \cd (\DDD \cd A) =0$ in addition to $\DDD \w (\DDD \w A)=0$, so that $\DDD^2$ is in general a grade-preserving operator. Duality relations also show that the differential forms codifferential $\delta A \equiv * \,d * A = \DDD \cd A$ is equivalent to divergence.

\section{Concluding remarks}

This article has attempted to lay an elementary foundation for studying geometric calculus on pseudo-Riemannian manifolds. The framework was seen to provide conceptually simple proofs of many concepts, such as obtaining a simple derivation of the Levi-Civita coefficients in arbitrary non-holonomic bases, demonstrating how the tensor covariant derivative follows from a chain rule, and clarifying the importance of gradient bases as the reciprocal to holonomic bases. We rederived these known concepts within a new constructive framework---an axiomatic formalism more similar to standard Riemannian geometry and general relativity than existing treatments. We hope this can serve both as a pedagogical introduction to the topic, and as a useful bridge for physicists to apply geometric calculus concepts in the context of relativity.


\subsection*{Acknowledgements}

I would like to thank Anthony Aguirre, Adam Reyes, and Ross Greenwood for useful discussions, and to thank David Constantine for first teaching me Riemannian geometry. This research was supported by the Foundational Questions Institute (FQXi.org), of which AA is Associate Director, and by the Faggin Presidential Chair Fund.

\clearpage

\appendix

\section{Geometric algebra review}
\label{sec:GA}

Here we provide a basic review of geometric algebra (GA), both as an introduction and to set notation. For more detailed reviews, there are a number of useful options: the books and survey by Macdonald~\cite{macdonald2017,macdonald11,macdonald12} give a clear and basic introduction, the book by Doran~and Lasenby~\cite{doran07} highlights many applications to physics, and the original monograph by Hestenes and~Sobczyk~\cite{hestenes87} provides a more complete theoretical development. These each use somewhat different conventions, but are still mutually intelligible. Conventions here are in line with Macdonald.

The basic objects of geometric algebra are multivectors, which can be visualized as (sums of) \mbox{$k$-dimensional} parallelepipeds sitting in $N$-dimensional space. These parallelepipeds range from dimension zero (scalars), and one (vectors), up to $N$, the dimension of the ambient space. In this way multivectors provide a natural way to represent oriented lengths, areas, and volumes in physical space. Geometric algebra is also known as Clifford algebra (Clifford himself used the former name~\cite{Clifford1878}), and it includes exterior algebra as a sub-algebra under the wedge product.
 
In more detail, a geometric algebra is a linear space which can be built on top of an inner product vector space, by adding in the additional operation of wedge multiplication. The dot (inner) product
$$a \cd b$$
of vectors is a scalar, while the wedge (outer) product 
$$a \w b$$
of two vectors is a new type of object called a $2$-vector. In GA scalars ($0$-vectors), vectors \mbox{($1$-vectors)}, $2$-vectors (also called bivectors), and higher order $k$-vectors (formed by additional wedging, and visualized as $k$-dimensional parallelepipeds) are all a part of the same linear space, which is called the space of \textit{multivectors}. The order $k$ of a $k$-vector is called its \textit{grade}, and $k$-vectors of different grade are linearly independent from one another. The wedge product is antisymmetric on vectors and associative on all multivectors (thus totally antisymmetric under vector swaps), so that \mbox{$k$-vectors $ a_1 \w \ldots \w a_k$} can be formed all the way up to the dimension~$N$ of the vector subspace. An essential aspect of GA is that $k$-vectors can be linearly independently added, so that an arbitrary multivector is written
\begin{equation}
    A = \sum_{k=0}^{N} \grade{A}_k \; ,
\end{equation} 
where $\grade{}_k$ is called the \textit{grade operator}, which picks out the $k$-vector part of a multivector. When different grades are added together they do not combine, and should be visualized as the formal sum of two separate parallelepipeds. The ability to add together different grades is leveraged to combine the dot and wedge products into a more fundamental operation called the \textit{geometric product}, which acts on vectors $a$ and $b$ by
\begin{equation}
    ab = a \cd b + a \w b, 
\end{equation} 
so that the geometric product $ab$ of two vectors is equal to a scalar plus a $2$-vector. More generally, the geometric product of arbitrary of multivectors is written
$$AB$$
and from a technical perspective is considered the fundamental operation of GA (not derivable from other operations, an elementary construction is given by \cite{Macdonald2002}). The geometric product is associative and distributes over addition, but is not commutative. The dot and wedge products for general multivectors, meanwhile, are defined in terms of the geometric product. For pure-grade multivectors $A_j$~(a $j$-vector) and $B_k$~(a $k$-vector) the dot and wedge are defined by
\begin{equation}
    \begin{array}{rcl}
    A_j \wcd B_k &=& \grade{A_j B_k}_{\,k-j} \; , \\[2pt]
    A_j \w   B_k &=& \grade{A_j B_k}_{\,k+j} \; ,
    \end{array}
\end{equation}
so that the dot (wedge) product acts as the maximally-grade-lowering (maximally-grade-raising) component of the geometric product. This definition is extended to arbitrary multivectors, which can be written as sums of pure-grade multivectors, by linearity. On vectors the dot and wedge amount to
\begin{equation}
   \begin{array}{rcl}
    a \wcd b &=& (ab + ba)/2  \; , \\[2pt]
    a \w b &=& (ab - ba)/2   \; ,
    \end{array}
\end{equation}
which could also be derived from the earlier expression for $ab$ using \mbox{$a \cd b = b \cd a$} and \mbox{$a \w b = -b \w a$}. There is another useful identity for vectors \cite{macdonald2017}, written
\begin{equation}
    \label{eqn:ext-fund-iden}
    aB = a \cd B + a \w B,
\end{equation}
but in general
\begin{equation}
    AB \neq A \cd B + A \w B \,.
\end{equation}

There is no need to introduce a basis to work in GA, since one typically works directly with vectors and multivectors, rather than with components. But introducing a basis is sometimes useful. One typically starts with a vector basis $\e_i$ which can be extended to a canonical multivector basis $\e_J$ by wedging the vector basis elements (more details below in the discussion of geometric manifolds). Every vector basis $\e_i$ has a unique reciprocal basis $\e^i$ such that
\begin{equation*}
    \e^i \cd \e_j = \dd^i_j \; .
\end{equation*}
Which one has the upper and which one has the lower index is arbitrary, this is simply a pair of two vector basis sets which are mutually reciprocal. Reciprocality is a property of basis sets, and there is no such thing as the reciprocal of an individual vector (although a vector does have reciprocal \textit{components} $a=a^i \, \e_i = a_i \, \e^i$ relative to a pair of reciprocal bases). Even when a basis is introduced, all expressions in GA are independent of the basis choice. For example a vector $a$ can be decomposed with respect to two different bases $\e_i$ and $\tilde{\e}_i$ in the same equation without confusion, 
\begin{equation*}
\begin{array}{rcrcr}
    a &=& (a \cd \e^i) \, \e_i &=& (a \cd \tilde{\e}^i) \, \tilde{\e}_i \; \, \\[6pt]
      &=& a^i \, \e_i &=& \tilde{a}^i \, \tilde{\e}_i \; ,
\end{array}
\end{equation*}
as the equation itself is independent of any basis choice. The concept of ``change of basis" is not particularly important in the GA formalism, but if one desires to do so, one writes, for example, 
\begin{equation*}
    a = a^i \, \e_i =  a^i \, (\e_i \cd \tilde{\e}^j) \; \tilde{\e}_j = \tilde{a}^j \, \tilde{\e}_j
\end{equation*}
to go from basis $\e_i$ to basis $\tilde{\e}_i$. Components are used, but they are never removed from the context of well-defined expressions, where they always explicitly multiply basis vectors. Equalities are always written as equalities of multivectors and not equalities of components.

The possibility of adding together multivectors of various grades, combined with the identification of the geometric product as the fundamental type of multiplication, makes GA an especially useful framework in physics, geometry, and calculus~\cite{macdonald2017,doran07}.

\section{Smooth tangent structure}
\label{sec:smooth}

To ensure a suitable smooth structure for $GT_p M$ it must be constructed in terms of a coordinate system $x^i$ with coordinate tangent basis $\e(x_i)$. In terms of this basis the metric coefficients \mbox{$\e(x_i) \cd \e(x_j) = g_{ij}$} must be smooth coefficients forming a symmetric invertible matrix (with smooth inverse coefficients $g^{ij}$ defined by $g^{ij}g_{jk}=\dd^i_k$) at each point.

The coordinate vector basis $\e(x_i)$ (denote $\e(x_i) \equiv \E_i$ for convenience) generates a canonical multivector basis at each point in the following standard way. Let $J$ be a (nonrepeating, unordered) set of valid indices for $\E_i$. Denote by $\E_J$ the wedge product 
\begin{equation}
    \E_J=\E_{j_0} \w \ldots \w \E_{j_n}
\end{equation}
with the product taken over all $j_i \in J$, such that the indices $j_i$ strictly increase towards the right (when $J$ is the empty set $\E_J=1$). The set of $\E_J$ for all index sets $J$ forms a basis for multivectors in $GT_p M$. Every multivector field $A$ can be uniquely expressed in the canonical coordinate basis by
\begin{equation}
    A = A^J \E_J
\end{equation}
where $A^J$ are scalar field coefficients. A multivector field is considered \textit{smooth} when the coefficients $A^J$ are smooth in the multivector coordinate basis. This implies that a vector field $a$ is smooth if and only if the scalar fields $(a \cd \E^i)$ are smooth.

Consider a collection of vector fields $\e_i$ which form a basis for $T_p M$ at each point. There exists at each point a reciprocal basis $\e^i$ such that $\e^i \cd \e_j = \dd^i_j$. These $\e_i$ are a \textit{smooth (vector) basis frame} (or \textit{smooth basis vector fields}) if the scalar fields $(\e_i \cd \E^j)$ and $(\e^i \cd \E_j)$ are all smooth relative to the coordinate basis. 

A smooth vector frame generates a multivector frame $\e_J=\e_{j_0} \w \ldots \w \e_{j_n}$ analogous to the canonical coordinate frame above, and every multivector field has a unique expression $B = B^J \e_J$ in this frame. Under the the given definitions, it can be shown that multivector fields are smooth if and only if their coefficients are smooth in the $\e_J$ basis for any smooth vector frame $\e_i$.


\section{Proof of Theorem \ref{thm:thm:grade-pres-dd} (MDD preserves grade)}
\label{sec:app:mdd-grade-preserving-proof}

Let $D$ be a directional derivative, and let $\E_i$ be an orthonormal basis with associated multivector basis $\E_J$. It suffices to show that $D_{\E_i} \E_J$ is grade-preserving for all basis vectors $\E_i$ and basis multivectors $\E_J$. The proof follows from a simple pattern, where potentially non-grade-preserving terms come in pairs which together vanish by metric compatibility, but there is no simple notation to show the proof algebraically. The pattern can quickly be observed by calculating $D_{\E_i} ( \E_1 \E_2 \E_3)$ and $D_{\E_i} ( \E_1 \E_2 \E_3 \E_4)$ in the orthonormal basis, in which examples the non-grade-preserving term is seen to explicitly vanish. The pattern generalizes as follows. 

Since the basis is orthonormal, $\E^i = \pm \E_i$ and $\E^i \E_i = 1$, and $j\neq k$ implies \mbox{$\E^j \E_k = -\E_k \E^j$}. An arbitrary multivector basis element of grade $n$ can be written 
$\E_J = \E_{j_1} \ldots \E_{j_n}$,
where
\mbox{$j_1 < \ldots < j_n$}.
Consider the first term $D_{\E_i} (\E_{j_1}) \ldots \E_{j_n}$ in the product rule expansion for $D_{\E_i} \E_J$. Decompose the derivative into components in the $\E^k$ basis to write this term as
$$\sum\nolimits_k (D_{\E_i} (\E_{j_1}) \cd \E_k)\E^k \ldots \E_{j_n}.$$
Since $\E^k=\pm \E_k$ is orthogonal to $\E_j$ whenever $k \neq j$, this term has grade $n$ for all terms in the sum where $k\notin (j_2, \ldots, j_n)$. Every term which does not preserve grade has $k=j_p$ for some $j_p\in (j_2, \ldots, j_n)$, in which case anticommuting $\E^k$ through basis vectors and utilizing $\E^k \E_{j_p}=1$ leaves a constant times the multivector term $\E_J$ with both $\E_1$ and $\E_{j_p}$ removed. Thus every term either has grade $n$ or grade $n-2$, and the terms with grade $n-2$ are proportional to $\E_J$ with two basis vectors removed from the product. This generalizes to every term in the product rule expansion, so to show that grade is preserved one must show that the total grade $n-2$ term in the product rule expansion vanishes.

Denote by $J(p,q)$ the set of indices $(j_1,\ldots,j_n)$ with $j_p$ and $j_q$ removed. The preceding paragraph implies that
$D_{\E_i} \E_J = \grade{D_{\E_i} \E_J}_n + \grade{D_{\E_i} \E_J}_{n-2}$
where
$\grade{D_{\E_i} \E_J}_{n-2} = \sum_{(p,q)} A^{J(p,q)} \E_{J(p,q)}$.
We now show that $A^{J(p,q)}=0$ using metric compatibility.

Consider 
$D_{\E_i} (\E_{j_1} \ldots \E_{j_p} \ldots \E_{j_q} \ldots \E_{j_n})$.
Every term in the product rule expansion which does not take a derivative of either $\E_{j_p}$ or $\E_{j_q}$ does not contribute to $\E_{J(p,q)}$, since it can only eliminate one or the other of the two factors. Two terms remain. The $p$ derivative term is 
$\sum\nolimits_k \E_{j_1} \ldots (D_{\E_i} (\E_{j_p}) \cd \E_k)\E^k \ldots \E_{j_q} \ldots \E_{j_n}.$ 
This contributes to $\E_{J(p,q)}$ only for $k=q$. The $q$ derivative term is analogous, so that the total $\E_{J(p,q)}$ term is
$$
\Big( \E_{j_1} \ldots (D_{\E_i} (\E_{j_p}) \cd \E_{J_q})\E^{j_q} \ldots \E_{j_q} \ldots \E_{j_n} \Big)
+
\Big( \E_{j_1}  \ldots \E_{j_p} \ldots (D_{\E_i} (\E_{j_q}) \cd \E_{j_p})\E^{J_p} \ldots \E_{j_n} \Big)
.$$
Commuting scalars and anticommuting basis vectors makes this term equal to
$$ 
(-1)^{(q-p-1)} 
\Big[D_{\E_i}(\E_{j_p}) \cd \E_{j_q} + D_{\E_i} (\E_{j_q}) \cd \E_{j_p} \Big]
 \E_{J(p,q)}.
$$
Thus by metric compatibility $A^{J(p,q)} = \pm D_{\E_i} (\E_{j_p}\cd \E_{j_q}) = 0$ since the basis is orthonormal. 

Since every term proportional to $\E_{J(p,q)}$ is zero, and every other term in the product rule expansion preserves grade, the derivative obeys $\grade{D_{\E_i} \grade{\E_J}_n}_n =D_{\E_i} \grade{\E_J}_n$. Extending this result to arbitrary multivectors by linearity and the properties of $D$ implies that the directional derivative preserves the grade of all multivectors.


\section{Proof of Theorem \ref{thm:thm:mdd-affine-bijection} (MDDs $\leftrightarrow$ metric-compatible connections)}
\label{sec:app:mdd-existence-proof}

We must show that for any metric-compatible affine connection $\Dres$, there exists a multivector directional derivative operator $D$ which restricts to $\Dres$ when acting on vector fields. We achieve this by directly constructing an operator $D$ that restricts to $\Dres$, then showing that $D$ is an MDD by the axioms of Definition \ref{thm:def:directional-deriv}.

$D$ is constructed by defining its action on multivectors in a particular basis. Throughout the proof, $\E_i$ refers to the specific orthonormal basis in terms of which $D$ is defined. Likewise, $\E_J$ is the canonical multivector basis constructed from $\E_i$.

\renewcommand{\myspace}{6pt}

\begin{gbox}
\begin{defn}
\label{thm:def:D-operator}
Let $\Dres$ be a metric-compatible affine connection. 
Let $\E_i$ be a set of smooth orthonormal basis vector fields with associated canonical multivector basis $\E_J$.  Every smooth multivector field can be written uniquely in the form $A = \sum_J A^J \E_J$, where $A^J$ are smooth scalar field coefficients. Define an operator \mbox{$D: T_p M \times MV\!F(M) \to GT_p M$} by \\[4pt]
\begin{tabular}{rl}
(i) & $D_{(\alpha \, a + \beta \, b)} A = \alpha \, D_a A + \beta \, D_b A,  $ 
\\[\myspace]
(ii) & $D_{\E_i}(1) = 0$,
\\[\myspace]
(iii) & $D_{\E_i}(\E_j) = \Dres_{\E_i}(\E_j)$,
\\[\myspace]
(iv) & $D_{a}(\sum_J A^J \E_J) = \sum_J \left( ( \pp_a A^J ) \, \E_J + A^J \, ( D_{a} \E_J) \right) ,  $ 
\\[\myspace]
(v) & For $E_J$ with grade greater than one, $D_{a}(\E_{J})$ defined by the product rule:
\\[2pt]
 & If $i < k$ for all $k \in K$, then $D_{a}(\E_{i} \E_{K}) = (D_{a}\E_{i}) \E_{K} + \E_{i} ( D_{a} \E_{K}) $.
\end{tabular} \\[6pt]
This suffices to uniquely define the operator $D$ on all smooth multivector fields, and $D$ restricts to $\Dres$ when acting on vector fields.
\end{defn}
\end{gbox}

This $D$ is a well-defined operator which is equal to the metric-compatible affine connection $\Dres$ when restricted to act on vector fields. If $D$ satisfies all the axioms of Definition \ref{thm:def:directional-deriv}, then $D$ is a multivector directional derivative. All the axioms except the product rule come easily.

\begin{gbox}
\begin{lem}
\label{thm:lem:gamma-op-props}
$D$ satisfies properties (i-iv) of Definition \ref{thm:def:directional-deriv}.
\end{lem}
\begin{pf}
(i) True by assumption. (ii) Let $\grade{A}_0 = \ph$. Then 
$D_a \grade{A}_0 = D_a (\ph 1) = (\pp_a \ph )1 + \ph (D_a 1) = \pp_a \ph =  \pp_a \grade{A}_0 $. (iii) Let $\grade{A}_1 = b^j \E_j$. Then  
$D_a \grade{A}_1 = 
D_a (b^j \E_j) = 
a^i (\pp_{E_i} b^j) \E_j + a^i b^j (D_{\E_i} \E_j) =
a^i (\pp_{E_i} b^j) \E_j + a^i b^j \gamma_{ijk} \E^k$, which is a vector. (iv) $D_a(A+B) = D_a(A^J \E_J + B^J \E_J) = D_a( (A^J+B^J)\E_J) = \pp_a(A^J+B^J) \E_J + (A^J+B^J) D_a \E_J = D_a A + D_a B$.
\end{pf}
\end{gbox}

Thus if $D$ obeys the product rule then it is a multivector directional derivative. 

The MDD axioms imply metric-compatibility on vectors. It follows that if $\Dres$ were not metric-compatible, then $D$ would have to violate the product rule. But $\Dres$ \textit{is} metric-compatible, and it turns out this is sufficient to ensure the product rule is obeyed by $D$. The method of proof is to show that the product rule holds as a special case for increasingly general strings of the orthonormal basis vectors $\E_i$ in terms of which the action of $D$ is defined.

First, a product rule is shown for pairs of basis vectors.

\begin{gbox}
\begin{lem}
On a pair of the defining basis vectors
\begin{equation*}
D_{\E_i} (\E_j \E_k) = ( D_{\E_i} \E_j ) \E_k  + \E_j (D_{\E_i} \E_k) \, .
\end{equation*}
\end{lem}
\begin{pf}
Define $\Dres_{\E_i} \E_j = \gamma_{ijk} \E^K$. 
Since $\Dres$ is metric-compatible and $\E_i$ orthonormal, $\gamma_{ijk}+\gamma_{ikj}=0$.
The proof follows from showing that 
$D_{\E_i} (\E_j \E_k) - ( D_{\E_i} \E_j ) \E_k  - \E_j (D_{\E_i} \E_k) = 0 $.
There are three cases. 
For $(j<k)$, the product rule is assured by definition of $D$. 
For $(j=k)$,  
$D_{\E_i} (\E_j \E_j) - ( D_{\E_i} \E_j ) \E_j  - \E_j (D_{\E_i} \E_j) = 
0 - 2 (D_{\E_i} \E_j) \cd \E_j =
-2 \gamma_{ijj} = 0$.
For $(j>k)$, one finds 
$D_{\E_i} (\E_j \E_k) - ( D_{\E_i} \E_j ) \E_k  - \E_j (D_{\E_i} \E_k) = 
-D_{\E_i} (\E_k \E_j) - ( D_{\E_i} \E_j ) \E_k  - \E_j (D_{\E_i} \E_k) =
- (D_{\E_i} \E_k) \E_j - \E_k (D_{\E_i} \E_j)  - ( D_{\E_i} \E_j ) \E_k  - \E_j (D_{\E_i} \E_k) =
  - 2 ( D_{\E_i} \E_j ) \cd \E_k   - 2 ( D_{\E_i} \E_k ) \cd \E_j  =
  -2 (\gamma_{ijk} + \gamma_{ikj}) = 0$.
\end{pf}
\end{gbox}

This can be extended from pairs of basis vectors to strings of basis vectors.

\renewcommand{\myspace}{4pt}
\begin{gbox}
\begin{lem} 

Say that an arbitrary string 
$\E_{j_1} \ldots \E_{j_n} $
of $n$ of the defining basis vectors is \textit{product-rule-separable} if it splits into derivatives of individual basis vectors in the usual way:
\begin{equation*}
D_{\E_i} (\E_{j_1} \ldots \E_{j_n}) = 
(D_{\E_i} \E_{j_1}) \E_{j_2} \ldots \E_{j_n} +
 \ldots +
\E_{j_1}  \ldots \E_{j_{n-1}}  (D_{\E_i} \E_{j_n}) .
\end{equation*}
Then product rule separability of basis strings is preserved by the basic manipulations of neighbor-swapping and pair insertion. That is, if $\epsilon = \E_{j_{0}} \ldots \E_{j_{m}}\E_{j_{m+1}} \ldots \E_{j_{n}}$ is product-rule-separable, then

\vspace*{\myspace}
$\,$(i) $\epsilon' = \E_{j_{0}} \ldots \E_{j_{m+1}} \E_{j_{m}} \ldots \E_{j_{n}}$ is product-rule-separable,

\vspace*{\myspace}
(ii) $\epsilon'' = \E_{j_{0}} \ldots \E_{j_{m}}  \E_k \E_k  \E_{j_{m+1}} \ldots \E_{j_{n}}$ is product-rule-separable, for any $\E_k$.
\end{lem}
\begin{pf}
(i) If $\E_{j_m} = \E_{j_{m+1}}$ then $\epsilon' =  \epsilon$ which is separable. Otherwise, $D_{\E_i} \epsilon' = D_{\E_i}(- \epsilon) = - D_{\E_i} \epsilon $. One can show that $ - D_{\E_i} \epsilon $ is equal to the product rule expansion for $\epsilon'$. Expand $D_{\E_i} \epsilon$ by the product rule. The product rule expansion for $\epsilon'$ is obtained by swapping $m \leftrightarrow m+1$ in each term. This yields a negative sign in each term, seen as follows. For each term with $\E_m \E_{m+1}$ outside the derivative, the identity $\E_j \E_k = - \E_k \E_j $ for $j\neq k$ suffices to give the negative sign. The two remaining terms are of the form 
$A [ (D_{\E_i}\E_m) \E_{m+1} + \E_m (D_{\E_i} \E_{m+1} ) ]B = 
A(D_{\E_i}(\E_m \E_{m+1}))B = 
A(D_{\E_i}(- \E_{m+1} \E_m))B = 
- A(D_{\E_i}(\E_{m+1} \E_m))B = 
-A [ (D_{\E_i}\E_{m+1}) \E_{m} + \E_{m+1} (D_{\E_i} \E_{m} ) ]B$,
using the product rule proved earlier for pairs of basis vectors, which yields the desired negative sign. Thus $D_{\E_i} \epsilon'$ is equal to its product rule expansion. 
(ii) $D_{\E_i} \epsilon'' = D_{\E_i}(\E_k^2 \epsilon) = \E_k^2 D_{\E_i} \epsilon $ since $\E_k^2$ is a constant scalar. Expand $E_k^2 D_{\E_i} \epsilon$ by the product rule. For every term with no derivatives of $\E_k$, the scalar $\E_k^2$ can be freely moved between $\E_m$ and $\E_{m+1}$. The remaining two terms are obtained by adding in a term of the form 
$ 0 =
A(D_{\E_i}(\E_k \E_k))B = 
A [ (D_{\E_i}\E_k) \E_{k} + \E_{k} (D_{\E_i} \E_{k} ) ]B
$
using the product rule for pairs of basis vectors. Thus $E_k^2 D_{\E_i} \epsilon$ is equal to the product rule expansion of $D_{\E_i} \epsilon''$.
\end{pf}
\end{gbox}

Using these operations, product rule separability can be shown for arbitrary basis strings, which implies the product rule for arbitrary multivectors.

\renewcommand{\myspace}{4pt}
\begin{gbox}
\begin{lem}
It follows from the previous lemma that

\vspace*{\myspace}
$\,\,$(i) Every string of the basis vectors $\E_i$ is product-rule-separable,

\vspace*{\myspace}
$\,$(ii) $D_{\E_i} ( \E_J \E_k ) = (D_{\E_i}  \E_J ) \E_k  +  \E_J ( D_{\E_i}  \E_k  )$ for all $\E_i, \E_J , \E_K$,

\vspace*{\myspace}
(iii) $D_{a} ( AB ) = (D_{a}  A ) B  +  A ( D_{a}  B  )$ for all $a,A,B$.
\end{lem}
\begin{pf}
(i) Every $E_J$ is product-rule-separable by definition of $D$. Every arbitrary string of basis vectors can be formed from some $E_J$ by a series of  neighbor-swaps and pair-insertions. These operations preserve product-rule-separability as shown above. 
(ii) $\E_J$, $\E_K$, and $\E_J \E_K$ are all strings of basis vectors, and thus can be expanded and recombined by their product rule expansions, as shown above. 
(iii) This follows from the definition of $D$ from direct calculation of $D_{a} (AB)$ by expansion in the canonical basis, using additivity from Lemma \ref{thm:lem:gamma-op-props} and the product rule for $\E_J \E_k$.
\end{pf}
\end{gbox}

Thus it has been shown that $D$ obeys the product rule and thus satisfies all axioms of Definition~\ref{thm:def:directional-deriv}. So for each metric compatible affine connection $\Dres$, $D$ is a multivector directional derivative restricting to it. This concludes the proof.

It is worth noting that, so far, it seems that making fundamental use of the neighbor-swapping and pair insertion operations in an orthonormal canonical basis is the only reasonable way to attack the above proof. Some time after reaching that conclusion, I found out that Alan Macdonald's simplest ``elementary construction'' of GA~\cite{Macdonald2002} (which long precedes this work) is based on exactly the same operations. This seems quite interesting, and points to the important role of these operations underlying the structure of GA.


\section{Gradient and linear transformations}
\label{sec:interp-grad}

Does the gradient operator encode all information about the directional derivative of a multivector field? No. This section explains why, by investigating similar operators in linear algebra. 

Let $\mathcal{A}$ be a geometric algebra, and let $f:\grade{\mathcal{A}}_1 \to \grade{\mathcal{A}}_1$
be a linear transformation of the vector subspace of $\mathcal{A}$, so that for vectors $a,b$ and scalars $\alpha, \beta$,
\begin{equation*}
     f(\alpha a + \beta b) = \alpha f(a) + \beta f(b) .
\end{equation*}
Let $\e_i$ be an arbitrary vector basis such that $\e_i \cd \e_j = g_{ij}$, with reciprocal basis $\e^i$ such that $\e^i \cd \e_j = \dd^i_j$ and $\e^i = g^{ij} \, \e_j$, where $g^{ij}$ is the matrix inverse to $g_{ij}$. By linearity $f$ is determined by its scalar components, defined by
$$ f(\e_i) = f_{ij} \; \e^j = {f_{i}}^{j} \; \e_j \; ,
\qquad \qquad
f(\e^i) = {f^{i}}_{j} \; \e^j = f^{ij} \; \e_j \; ,$$
or equivalently by its vector components, defined by
$$ f(\e_i)=f_i \; , 
\qquad \qquad 
f(\e^i)=f^i \; . $$
It follows from linearity that $g_{ij}$ and its inverse can raise and lower each index of the components. 

Given $f$, define a multivector $\f$ by (with summation convention)
$$ 
\begin{array}{rcl}
    \f &=& \e^i \, f(\e_i) \\[4pt]
       &=& \e^i \cd f_i \; + \; \e^i \w f_i \\[4pt]
       &=& \tr(f) \; + \; \rot(f)
\end{array}
$$
where the \textit{trace} $\tr(f)$ and \textit{rotation} $\rot(f)$ of $f$ are defined by
$$
\tr(f) = \e^i \cd f(\e_i) \; ,
\qquad \qquad
\rot(f) = \e^i \w f(\e_i) \; .
$$
The trace is a scalar and the rotation is a bivector. In components,
$$
\f = f^{ij} \, (g_{ij} + \e_i \w \e_j) \; ,
\qquad \qquad
\tr(f) = f^{ij} \, g_{ij} \; ,
\qquad \qquad
\rot(f) = f^{ij} \; \e_i \w \e_j  \; .
$$
The trace, rotation, and $\f$ can all be checked to be independent of what basis they are defined in. Note that the trace depends only on $(f^{ij}+f^{ji})$ and the rotation depends only on $(f^{ij}-f^{ji})$. Geometrically, the trace represents the average amount by which input vectors are scaled parallel to themselves by $f$, while the rotation encodes how much input vectors contribute to a component purely perpendicular to themselves.

The trace and rotation are best understood by noticing that every linear transformation admits the unique decomposition
$$
f = \left( \tfrac{\tr(f)}{n} \right) \id  + f^{+} + f^{-}
$$
where $\id$ is the identity operator, $n$ is the dimension of the space of vectors, and $f^{\pm}$ are traceless symmetric ($+$) and antisymmetric ($-$) linear transformations satisfying $\left(f^{\pm}(a) \cd b \right) = \pm \left( a \cd f^{\pm}(b) \right)$ and $\tr(f^{\pm})=0$. This can be called the \textit{TSA (trace-symmetric-antisymmetric) decomposition}.

Clearly $\tr(f)$ encodes the strength of the identity component in $f$. More surprising, perhaps, is that the bivector $\rot(f)$ encodes the entirety of the antisymmetric part. In fact, for all $a$,
$$ a \cd \rot(f) = f^-(a) \; .$$
It follows that in general
$$ f(a) = \left( \tfrac{\tr(f)}{n} \right) a + a \cd \rot(f)  + f^{+}(a) \; , $$
and therefore that $\f = \tr(f)+\rot(f)$ encodes the $\id$ and $f^-$ parts of $f$ but not the $f^+$ part. The remaining part $f^{+}$ also admits a simple interpretation, since it can be shown that every symmetric transformation has a basis of eigenvectors and obeys the spectral theorem~\cite{macdonald11}. Thus $f^{+}$ is specified by a basis  $v_i$ of eigenvectors and a set of eigenvalues $\lambda_i$, in terms of which 
$$
f^+(a) = \sum_i \lambda_i \, (a \cd v^i) \, v_i \; .
$$ \\[-10pt]
Since $f^+$ is traceless, the eigenvalues must sum to zero ($\sum_i \lambda_i = 0$), which expresses in some heuristic sense that $f^{+}$ is ``conservative".

To summarize, a linear transformation $f$ is the sum of three basic operations. A trace term, characterized by a scalar $\tr(f)$, multiplies inputs by a constant scale~factor. A rotation term, characterized by a bivector $\rot(f)$, returns a vector orthogonal to the input. And a traceless-symmetric term, characterized by an eigenbasis $v_i$ with eigenvalues $\lambda_i$ summing to zero, scales inputs along its characteristic directions. The multivector 
$$\f=\e^i \, f(\e_i) = \tr(f)+\rot(f)$$ 
encodes the trace and rotation terms, but not the traceless-symmetric term. 

For a fixed vector field $b$, the directional derivative $f(a)=D_{a} b$ is a linear transformation from vectors to vectors. The gradient $\D b = \e^i D_{\e_i} b$ is exactly analogous to $\f$ for this transformation. The gradient of a vector field therefore encodes the trace (divergence) and rotation (curl) of the directional derivative, but ignores the traceless-symmetric part. This is specific to the gradient of vector fields---for scalar fields the gradient encodes all information since $D_a \ph = a \cd \D \ph$, while for higher grade fields the picture is more complicated.


\bibliographystyle{unsrt}
\bibliography{geom.bib}

\begin{thebibliography}{10}

\bibitem{Clifford1878}
Professor Clifford.
\newblock Applications of grassmann's extensive algebra.
\newblock {\em American Journal of Mathematics}, 1(4):350--358, 1878.

\bibitem{hestenes87}
D.~Hestenes and Garret Sobczyk.
\newblock {\em Clifford Algebra to Geometric Calculus: A Unified Language for
  Mathematics and Physics (Fundamental Theories of Physics)}.
\newblock Springer, 1987.

\bibitem{doran07}
Chris Doran and Anthony Lasenby.
\newblock {\em Geometric Algebra for Physicists}.
\newblock Cambridge University Press, 2007.

\bibitem{macdonald2017}
Alan Macdonald.
\newblock A survey of geometric algebra and geometric calculus.
\newblock {\em Advances in Applied Clifford Algebras}, 27(1):853--891, Mar
  2017.

\bibitem{lawson1989}
H.~Blaine Lawson, Jr. and Marie-Louise Michelsohn.
\newblock {\em Spin geometry}, volume~38 of {\em Princeton Mathematical
  Series}.
\newblock Princeton University Press, Princeton, NJ, 1989.

\bibitem{crumeyrolle1990}
Albert Crumeyrolle.
\newblock {\em Orthogonal and symplectic {C}lifford algebras}, volume~57 of
  {\em Mathematics and its Applications}.
\newblock Kluwer Academic Publishers Group, Dordrecht, 1990.

\bibitem{rodrigues2007}
Waldyr~Alves Rodrigues, Jr. and Edmundo Capelas~de Oliveira.
\newblock {\em The many faces of {M}axwell, {D}irac and {E}instein equations},
  volume 722 of {\em Lecture Notes in Physics}.
\newblock Springer, Berlin, 2007.

\bibitem{hestenes1984}
David Hestenes and Garret Sobczyk.
\newblock {\em Clifford algebra to geometric calculus}.
\newblock Fundamental Theories of Physics. D. Reidel Publishing Co., Dordrecht,
  1984.

\bibitem{doran2003}
Chris Doran and Anthony Lasenby.
\newblock {\em Geometric algebra for physicists}.
\newblock Cambridge University Press, Cambridge, 2003.

\bibitem{lee02}
John~M. Lee.
\newblock {\em Introduction to Smooth Manifolds (Graduate Texts in
  Mathematics)}.
\newblock Springer, 2002.

\bibitem{lee97}
John~M. Lee.
\newblock {\em Riemannian Manifolds: An Introduction to Curvature (Graduate
  Texts in Mathematics)}.
\newblock Springer, 1997.

\bibitem{Hestenes2015}
David Hestenes.
\newblock {\em Space-time algebra}.
\newblock Birkh\"{a}user/Springer, Cham, second edition, 2015.

\bibitem{macdonald11}
Alan Macdonald.
\newblock {\em Linear and Geometric Algebra}.
\newblock CreateSpace Independent Publishing Platform, 2011.

\bibitem{macdonald12}
Alan Macdonald.
\newblock {\em Vector and Geometric Calculus}.
\newblock CreateSpace Independent Publishing Platform, 2012.

\bibitem{Schutz1980}
Bernard {Schutz}.
\newblock {\em {Geometrical Methods of Mathematical Physics}}.
\newblock {Cambridge University Press}, 1980.

\bibitem{Macdonald2002}
Alan Macdonald.
\newblock An elementary construction of the geometric algebra.
\newblock {\em Advances in Applied Clifford Algebras}, 12(1):1--6, Jun 2002.

\end{thebibliography}

\end{document}